\documentclass[english,reqno,11pt]{amsart}

\usepackage[letterpaper,margin=2cm]{geometry}
\usepackage{amssymb,amsmath,amsthm,esint,braket,mathtools,color}

\usepackage{hyperref}

\usepackage[abbrev,msc-links]{amsrefs}
\usepackage{mathscinet}

\usepackage[capitalise]{cleveref}

\usepackage[T1]{fontenc}
\usepackage{lmodern}
\providecommand{\abs}[1]{\left\vert #1 \right\vert}
\providecommand{\norm}[1]{\left\Vert #1 \right\Vert}

\providecommand{\pt}[1]{\left( #1 \right)}
\providecommand{\spt}[1]{\left[ #1 \right]}

\providecommand{\RR}{\mathbb R}

\providecommand{\D}{\,\mathrm d}

\renewcommand{\epsilon}{\varepsilon}

\renewcommand{\qed}{\hfill$\blacksquare$\vspace{1ex}}
	
\theoremstyle{plain}
\newtheorem{thm}{Theorem}
\newtheorem{prop}{Proposition}
\newtheorem{lem}{Lemma}

\newtheorem{thmx}{Theorem}

\theoremstyle{remark}
\newtheorem{remark}{Remark}

\Crefname{cor}{Corollary}{Corollaries}
\Crefname{lem}{Lemma}{Lemmas}
\Crefname{thm}{Theorem}{Theorems}
\Crefname{thmx}{Theorem}{Theorems}
\Crefname{ineq}{Ineq.}{Ineqs.}
\creflabelformat{ineq}{(#2#1#3)}
	
\hypersetup{
	unicode=true,          % non-Latin characters in Acrobat’s bookmarks
	pdftitle={Simple proofs of a {G}agliardo--{N}irenberg inequality and their extensions},    % title
	pdfauthor={Hern{\'a}n Castro},     % author
	pdfcreator={Hern{\'a}n Castro},   % creator of the document
}
	
\begin{document}

\title{Simple proofs of a {G}agliardo--{N}irenberg inequality and their extensions}
\author{Hern{\'a}n Castro}
\date{\today}
\email{hcastro@utalca.cl}
\address{Instituto de Matem{\'a}ticas, Universidad de Talca, Casilla 747, Talca, Chile}

\begin{abstract}
We revisit the one-dimensional Gagliardo--Nirenberg interpolation inequality
\[
\Vert{u'}\Vert_{L^p(\RR)}\leq C\Vert{u}\Vert_{L^q(\RR)}^{1-a}\Vert{u''}\Vert_{L^r(\RR)}^{a},\qquad u\in C^\infty_c(\RR),\quad\tfrac12\leq a\leq1,
\]
and recall how its general form follows from the case \(a=\frac12\) by interpolation with an endpoint estimate. For \(a=\frac12\) and \(p=q=r\), we give four alternative proofs, based on integration by parts (two proofs), on an ODE, and on the \(L^p\) theory of singular integrals. The two proofs by integration by parts also cover most of the other admissible exponents. In \(\RR^N\), the integration by parts arguments give the inequality with explicit constants and, for \(p\leq2\), with the Laplacian in place of the Hessian and a constant independent of the dimension. For weighted measures \(w(x)\D x\), a pointwise version of the ODE argument gives the pointwise inequality \(|\nabla u|^2\leq C_N\,Mu\,M(\Delta u)\), where \(M\) is the Hardy--Littlewood maximal function, and hence the inequality for Muckenhoupt weights, with the Laplacian and an explicit constant. The integration by parts arguments give it for homogeneous weights of large degree and for superharmonic weights. These results apply, in particular, to power weights and to monomial weights.
\end{abstract}


\maketitle

\section{Introduction}

In this article we revisit the following classical result, due to Gagliardo \cite{Gag59} and Nirenberg \cite{Nir59}. It is the case of one variable, and of the first and second derivatives, of the Gagliardo--Nirenberg interpolation inequalities (see the Theorem in \cite{Nir59}*{p. 125}).
\begin{thm}[Gagliardo--Nirenberg interpolation inequality]\label{gn-thm}
Let \(1\leq q,r\leq\infty\) and \(\frac12\leq a\leq1\), and assume that the number \(p\) given by
\begin{equation}\label{gen-GN2}
\frac{1}{p}=1+a\pt{\frac1r-2}+\frac{1-a}q
\end{equation}
satisfies \(\frac1p\geq0\). Then \(1\leq p\leq\infty\), and there exists a constant \(C>0\), depending only on \(a\), \(q\) and \(r\), such that for every \(u\in C^\infty_c(\RR)\)
\begin{equation}\label{gen-GN}
\norm{u'}_{L^p(\RR)}\leq C\norm{u}_{L^q(\RR)}^{1-a}\norm{u''}_{L^r(\RR)}^a.
\end{equation}
\end{thm}
When \(a=\frac12\), \eqref{gen-GN2} reduces to \(\frac2p=\frac1q+\frac1r\), and \eqref{gen-GN} reads
\begin{equation}\label{GN-1D}
\norm{u'}_{L^p(\RR)}\leq C\norm{u}_{L^q(\RR)}^{\frac12}\norm{u''}_{L^r(\RR)}^{\frac12},
\end{equation}
where, by Nirenberg's proof, \(C\) can be taken independent of \(p\), \(q\) and \(r\) \cite{Nir59}*{p. 129}.

The Gagliardo--Nirenberg interpolation inequalities play a central role in analysis and partial differential equations, and they have been studied extensively since the works \cite{Gag59,Nir59}. It is not our intention to give a complete survey of these inequalities, but we refer the interested reader to the work of Brezis and Mironescu \cite{BrMi2018}, where the full range of validity of the Gagliardo--Nirenberg inequalities in fractional Sobolev spaces is described, and to \cite{FiFoRoSo2021} for a detailed version of Nirenberg's proof together with some historical remarks. Some particular cases of \cref{gn-thm} were known well before \cite{Gag59,Nir59}. For instance, if \(a=\frac12\) and \(p=q=r=\infty\), then \eqref{gen-GN} is the classical inequality of Landau \cite{Lan1913}
\begin{equation}\label{landau}
\norm{u'}_{L^\infty(\RR)}^2\leq 2\norm{u}_{L^\infty(\RR)}\norm{u''}_{L^\infty(\RR)},
\end{equation}
where the constant \(2\) is optimal. Kolmogorov \cite{Kol1939} later extended this inequality to higher-order derivatives. Similarly, if \(a=\frac12\) and \(p=q=r=2\), the analogue of \eqref{gen-GN} on the half-line \((0,\infty)\) is an inequality of Hardy and Littlewood (see \cite{HL52}*{Section 7.8}). We refer to the monograph of Kwong and Zettl \cite{KwZe1992} for many results of this kind, both on the real line and on the half-line. At the other end, \(a=1\) forces \(r=1\) and \(p=\infty\), and \eqref{gen-GN} is the one-dimensional Sobolev inequality for \(u'\) (see \cref{rem:general-a}).

\Cref{gn-thm} follows from its case \(a=\frac12\), that is, from \eqref{GN-1D}, by H{\"o}lder interpolation with an endpoint estimate. In dimension one, the endpoint estimate can be obtained with an ODE argument, which gives \eqref{gen-GN} directly when \(p\geq\max\{q,r\}\) (see \cref{sec:general-a}). In \(\RR^N\), with \(N\geq2\) and \(r<N\), it is Sobolev's inequality (see \cref{lem:general-a} and \cref{rem:general-ND}). Nirenberg \cite{Nir59}*{p. 126} also reduces the general case to the extreme values of \(a\) by interpolation. For this reason, in the rest of this work we focus on the case \(a=\frac12\), that is, on \eqref{GN-1D}, and we state our results for that case.

As one can see in \cite{Nir59}*{pp. 129--131}, Nirenberg's proof of \eqref{GN-1D} is very elementary in its ingredients; however, it is a bit tricky, as it relies on a clever way of partitioning the real line (see also \cite{GT01}*{Theorem 7.27} for another proof in the same spirit). As is often the case, one of the advantages of finding alternative proofs of the same result is that each proof provides different insights into the result itself. This is the leading goal of this work: to provide alternative proofs of \eqref{GN-1D}, using techniques that differ from the one used by Nirenberg, and to use them to understand \eqref{GN-1D} and its extensions better. More precisely, in \cref{sec:p=q=r} we give four proofs of \eqref{GN-1D} when \(p=q=r\), each of them based on a different classical tool. We present them according to the tools they require, from the most elementary to the most advanced:
\begin{itemize}
\item Proof \#1 is an integration by parts against \(\abs{u'}^{p-2}u'\), the one-dimensional \(p\)-Laplacian of \(u\), and it gives a stronger inequality when \(p\geq2\), which is well known (see \cite{KaPe2012}). When \(1<p<2\) it is combined with a duality argument.
\item Proof \#2 is an integration by parts against \(\abs{u}^{p-2}u\). This is the classical computation behind the \(L^p\)-dissipativity of the Laplacian. Its justification when \(1<p<2\) requires some care, because of the singularity of \(\abs{u}^{p-2}\) at the zeros of \(u\) (see \cite{MeSp2008}). The weighted estimate that it gives is due to Ka{\l}amajska and Peszek \cite{KaPe2012}.
\item Proof \#3 writes \(u'\) in terms of \(-u''+u\) by means of the fundamental solution of the ODE \(-u''+u=f\), which is the one-dimensional Bessel potential (see \cite{Ste70}*{Chapter V}). It concludes with Young's inequality and a scaling argument, and it covers every \(1\leq p\leq\infty\).
\item Proof \#4 relies on the \(L^p\) theory of singular integrals of Calder{\'o}n and Zygmund \cite{CaZy1952} (see also \cite{Ste70}), applied to the function \(u(x)u(y)\) of two variables. It is very short, but it requires these more advanced tools.
\end{itemize}
Let us stress that the ingredients of Proofs \#1 and \#2 are classical (see the references above), and that Proofs \#3 and \#4 rely on standard tools. Our purpose is to show that each of them gives a complete proof of \eqref{GN-1D}, with explicit constants in the first three cases, and that they extend in different directions.
Proofs \#1 and \#2 also give \eqref{GN-1D} for most of the exponents \(p,q,r\) that are not all equal. The remaining ones, which have \(r=1\), are covered by Nirenberg's result \cite{Nir59}, and also by a recent result of Le{\'s}nik, Roskovec and Soudsk{\'y} \cite{LeRoSo2025} (see \cref{prop:r=1}).

The first natural generalization of \eqref{GN-1D} is to consider functions of several variables. If one replaces \(u'\) by \(\nabla u\) and \(u''\) by \(D^2u\), inequality \eqref{GN-1D} remains valid in \(\RR^N\). This particular case of the Gagliardo--Nirenberg inequalities in several variables is due to Gagliardo \cite{Gag59} and Nirenberg \cite{Nir59}*{Inequality (2.5), p. 129} (see \cite{FiFoRoSo2021}*{Section 2} for the history), and it follows from the one-dimensional inequality by integrating over lines (\cref{gn-thm-ND} below). When \(1<r<\infty\), the Calder{\'o}n--Zygmund estimate \(\norm{D^2u}_{L^r(\RR^N)}\leq C\norm{\Delta u}_{L^r(\RR^N)}\) also gives it with \(\Delta u\) in place of \(D^2u\). Since Proofs \#1 and \#2 only use integration by parts and H{\"o}lder's inequality, they can be carried out directly in \(\RR^N\), and they give these inequalities with simple explicit constants.
\begin{thmx}\label{thmA}
Let \(N\geq1\), and let \(1<p<\infty\) and \(1\leq q,r\leq\infty\) satisfy \(\frac2p=\frac1q+\frac1r\). Then for every \(u\in C^\infty_c(\RR^N)\) the following holds:
\begin{enumerate}
\item[(a)] if \(p\geq2\), then
\[\norm{\nabla u}_{L^p(\RR^N)}^2\leq\sqrt{N-1+(p-1)^2}\,\norm{u}_{L^q(\RR^N)}\norm{D^2u}_{L^r(\RR^N)};\]
\item[(b)] if \(p\leq2\) and \(r>1\), then
\[\norm{\nabla u}_{L^p(\RR^N)}^2\leq\frac{r'}{q}\norm{u}_{L^q(\RR^N)}\norm{\Delta u}_{L^r(\RR^N)},\qquad\text{where }\frac1r+\frac1{r'}=1.\]
\end{enumerate}
\end{thmx}
Here \(\abs{D^2u}\) denotes the Frobenius norm of the Hessian, \(\abs{D^2u}^2=\sum_{j,k=1}^N\abs{\partial_j\partial_ku}^2\). For \(N=1\) the constants are those of Proofs \#1 and \#2. For \(p\geq2\), the explicit constant \(\sqrt N+p-2\) also follows from an identity of G{\"u}neysu and Pigola \cite{GuPi2019}*{Theorem 1}, and part (a) improves it to \(\sqrt{N-1+(p-1)^2}\). In part (b), the constant \(\frac{r'}q\) does not depend on the dimension, and the Laplacian is obtained directly, without the Calder{\'o}n--Zygmund theory. For \(p=q=r\), dimension-free constants can also be derived from the dimension-free bounds for the Riesz transforms \cite{Ste1983}, but not by an elementary argument. \Cref{thmA} is \cref{gn-thm-RN} in \cref{sec:dimN}.

Using the Caffarelli--Kohn--Nirenberg inequalities \cite{CKN1984} as an inspiration, one can also ask for the validity of the \(N\)-dimensional inequality when the Lebesgue measure is replaced by a weighted measure \(w(x)\D x\). Writing \(L^t(w):=L^t(\RR^N,w\D x)\), one would like to have
\begin{equation}\label{GN-intro-w}
\norm{\nabla u}_{L^p(w)}^2\leq C\norm{u}_{L^q(w)}\norm{D^2u}_{L^r(w)},
\end{equation}
where \(\frac2p=\frac1q+\frac1r\), for instance for the power weights \(w(x)=\abs x^\alpha\). For power weights, inequalities like \eqref{GN-intro-w} are contained in the very general family of weighted interpolation inequalities obtained by C.-S. Lin \cite{Lin1986}, who extended the Caffarelli--Kohn--Nirenberg family to higher-order derivatives (see also \cite{DuSi2022} for a recent generalization). For more general weights, a useful tool is the family of pointwise multiplicative inequalities of Ka{\l}amajska \cite{Kal1994} and of Maz'ya and Shaposhnikova \cite{MaSh1999}, which bound \(\abs{\nabla u(x)}^2\) by the product of the Hardy--Littlewood maximal functions of \(u\) and of \(D^2u\) at \(x\). Combined with Muckenhoupt's theorem on the boundedness of the maximal function on weighted \(L^p\) spaces \cite{Muck1972}, these pointwise inequalities give \eqref{GN-intro-w} for Muckenhoupt weights (see \cite{Kal1994,KaPi2006}). More recently, Le{\'s}nik, Roskovec and Soudsk{\'y} \cite{LeRoSo2025} replaced the maximal function by sparse averaging operators, which gives the Gagliardo--Nirenberg inequality in every rearrangement-invariant Banach function space. Similarly, the weighted Calder{\'o}n--Zygmund theory of Coifman and Fefferman \cite{CoiFe1974} allows one to adapt Proof \#4 to Muckenhoupt weights. Other weights that have attracted attention recently are the monomial weights \(x^A=\abs{x_1}^{A_1}\cdots\abs{x_N}^{A_N}\) (see for instance \cite{CR-O2013-2,Cas2016-2,Cas2021}) and, more generally, homogeneous weights on convex cones \cite{CRS2016,CFR2020}. Many of these weights are not Muckenhoupt weights, and the inequalities associated with them are usually obtained with techniques that are specific to the weight under consideration.

With the above in consideration, in \cref{sec:dimN,sec:weights} we analyze which of the proofs of \eqref{GN-1D} can be extended to these settings.

In \cref{sec:weights} we consider weighted measures. A pointwise version of Proof \#3 gives the following result for Muckenhoupt weights.
\begin{thmx}\label{thmB}
Let \(\beta_N:=\sqrt\pi\,\Gamma\pt{\frac{N+1}2}/\Gamma\pt{\frac N2}\), let \(1<q,r\leq\infty\), let \(1<p\leq\infty\) be given by \(\frac2p=\frac1q+\frac1r\), and let \(w\in A_{\min\{q,r\}}\) (no condition on \(w\) is needed if \(q=r=\infty\)). Then for every \(u\in C^\infty_c(\RR^N)\)
\[\norm{\nabla u}_{L^p(w)}^2\leq4\beta_N^2\norm{M}_{L^q(w)}\norm{M}_{L^r(w)}\norm{u}_{L^q(w)}\norm{\Delta u}_{L^r(w)}.\]
\end{thmx}
Here \(M\) is the Hardy--Littlewood maximal function, \(A_s\) is the Muckenhoupt class, and \(\norm{M}_{L^s(w)}\) is the norm of \(M\) on \(L^s(w)\), with \(\norm{M}_{L^\infty(w)}:=1\). For \(1<r<\infty\), an inequality of this form, with an unspecified constant, also follows from the weighted Gagliardo--Nirenberg inequalities of \cite{Kal1994,KaPi2006} and the weighted Calder{\'o}n--Zygmund estimate of \cite{CoiFe1974}; the point of \cref{thmB} is the constant. Its proof relies on the pointwise inequality
\[\abs{\nabla u(x)}^2\leq4\beta_N^2\,Mu(x)\,M(\Delta u)(x),\qquad x\in\RR^N,\]
(\cref{lem:pointwise}), combined with Muckenhoupt's theorem. Pointwise inequalities of this kind, with \(D^2u\) in place of \(\Delta u\), are due to Ka{\l}amajska \cite{Kal1994}*{Theorem 3} and to Maz'ya and Shaposhnikova \cite{MaSh1999}, who also obtained explicit constants. In dimension one, where \(\Delta u=u''\) and \(\beta_1=1\), \cref{lem:pointwise} gives the constant \(4\), while Maz'ya and Shaposhnikova obtained the sharp constant \(\frac83\) for a smaller maximal operator \cite{MaSh2002}*{inequality (16)} (see \cref{rem:Ap}). For \(N\geq2\), we have not found the pointwise inequality with the Laplacian in the literature; it is not covered by the framework of \cite{Kal1994}, and for \(N\geq3\) Maz'ya and Shaposhnikova \cite{MaSh2000} prove an inequality of this form with \((-\Delta)^{1/2}u\) in place of \(\nabla u\) (see \cref{rem:Ap}). \Cref{thmB} is contained in \cref{gn-thm-Ap} in \cref{sec:weights}, which also covers \(p=q=r=1\) and \(w\in A_1\).

Beyond Muckenhoupt weights, Proofs \#1 and \#2 give \eqref{GN-intro-w} for homogeneous weights when the degree of the weight is large enough (\cref{gn-thm-hom}). Proof \#2 also gives \eqref{GN-intro-w} with \(\Delta u\) in place of \(D^2u\) for superharmonic weights, when \(p\leq2\), with the same constant as in dimension one, which does not depend on the weight (\cref{gn-thm-superh}).

Typical examples are the power weights \(\abs x^\alpha\), which belong to \(A_s\) if and only if \(-N<\alpha<N(s-1)\), and more generally the homogeneous weights of degree \(\alpha\) that are positive and continuous on the unit sphere, which are comparable to \(\abs x^\alpha\). For them, \cref{thmB} applies when \(-N<\alpha<N(\min\{q,r\}-1)\), and together with \cref{gn-thm-hom} it covers every \(\alpha>-N\) when \(N\geq2\) and \(r<Nq\) (see \cref{rem:power}, where we also compare with \cite{Lin1986,DuSi2022}). \Cref{thmB} also applies to the monomial weights \(\abs{x_1}^{A_1}\cdots\abs{x_N}^{A_N}\) with \(-1<A_i<\min\{q,r\}-1\), since products of one-dimensional \(A_s\) weights are \(A_s\) weights (\cref{rem:Ap}). These weights degenerate on the unit sphere; for them, a slicing argument reduces the inequality with \(D^2u\) to dimension one, where it holds for every \(A_i>-1\) if \(1<q,r<\infty\) (\cref{rem:weights-other}).

The rest of this paper is divided as follows. In \cref{sec:p=q=r} we provide four different proofs of \eqref{GN-1D} in the case \(p=q=r\). We also show that the first two of them give \eqref{GN-1D} for most of the remaining exponents, we recall in \cref{prop:r=1} the known result for the remaining ones, and in \cref{sec:general-a} we prove \cref{gn-thm}. In \cref{sec:dimN} we consider functions of several variables, and in \cref{sec:weights} weighted measures.

\medskip
\noindent\textbf{AI usage disclosure statement.}
Artificial intelligence tools were used to assist with double-checking routine computations, with the literature search, and to improve the clarity and presentation of the manuscript. The mathematical results and arguments are the author's own, computations and references were independently verified by the author, who takes full responsibility for the content of the manuscript.

\section{Four proofs of the case \texorpdfstring{\(a=\frac12\)}{a=1/2}}\label{sec:p=q=r}

In this section we give four proofs of \eqref{GN-1D} when \(p=q=r\), and in \cref{sec:general-a} we prove \cref{gn-thm}. We point out that from now on, we will always assume that \(u\not\equiv 0\), and, unless otherwise stated, that \(u\in C^\infty_c(\RR)\) in this section and \(u\in C^\infty_c(\RR^N)\) in \cref{sec:dimN,sec:weights}. Also, unless otherwise stated, \(C>0\) will denote a constant that depends only on the parameters of the problem and that may change from line to line.

\subsection{Proof \#1 -- Simple integration by parts}\label{sec:proof1}

In this subsection we prove \eqref{GN-1D} for \(p=q=r\) and \(1<p<\infty\).

\begin{proof}
When \(2\leq p<\infty\), one can prove the following stronger inequality by just integrating by parts and using H{\"o}lder's inequality: there exists a universal constant \(C>0\) such that 
\begin{equation}\label{strong-esti}
\int_\RR \abs{u'(x)}^p\D x\leq C\int_{\RR}\abs{u(x)u''(x)}^{\frac p2}\D x.
\end{equation}
Indeed, observe that if \(p\geq 2\) then \(\abs{u'}^{p-2}u'\in C^1_c(\RR)\) with \((\abs{u'}^{p-2}u')'=(p-1)\abs{u'}^{p-2}u''\), and therefore one can integrate by parts and use H{\"o}lder's inequality to obtain
\begin{align*}
\int_{\RR} \abs{u'(x)}^p\D x&=\int_{\RR} \abs{u'(x)}^{p-2}u'(x)u'(x)\D x\\
&=-\int_{\RR} \pt{\abs{u'(x)}^{p-2}u'(x)}'u(x)\D x\\
&=-(p-1)\int_{\RR}\abs{u'(x)}^{p-2}u''(x)u(x)\D x\\
&\leq (p-1)\pt{\int_{\RR}\abs{u'(x)}^p\D x}^{1-\frac2p}\pt{\int_{\RR}\abs{u(x)u''(x)}^{\frac p2}\D x}^{\frac 2p},
\end{align*}
which gives \eqref{strong-esti}, and after using H{\"o}lder's inequality once more we obtain \eqref{GN-1D}. Inequality \eqref{strong-esti}, with the constant \(C=(p-1)^{\frac p2}\) given by this computation, is well known (see \cite{KaPe2012}*{Lemma 1.1}).

As we show in \cref{rem:strong-fails} below, \eqref{strong-esti} is false for every \(1\leq p<2\). Nonetheless, one can use \eqref{strong-esti} together with a duality argument to prove \eqref{GN-1D} when \(1<p<2\) and \(p=q=r\). To do this, we need the following
\begin{lem}\label{lp-norm-equiv}
For all \(1<p<\infty\) and all \(v\in L^p(\RR)\), let
\begin{equation*} %\label{lp-norm-equiv-def}
\abs{v}_{L^p(\RR)}:=\sup_{\substack{\varphi\in W^{2,p'}(\RR)\\ \varphi\not\equiv 0}}\frac{\displaystyle\int_\RR v(x)\spt{\varphi(x)-\varphi''(x)}\D x}{\norm{\varphi}_{L^{p'}(\RR)}+\norm{\varphi''}_{L^{p'}(\RR)}},
\end{equation*}
where \(\frac1p+\frac1{p'}=1\). Then \(\abs{v}_{L^p(\RR)}\) defines an equivalent norm on \(L^p(\RR)\).
\end{lem} 

If we assume for a second that \cref{lp-norm-equiv} is true (we will prove this lemma in the appendix), we readily obtain a proof of \eqref{GN-1D} in the case \(1<p<2\). Indeed, apply \cref{lp-norm-equiv} to \(v=u'\), and integrate by parts to obtain
\begin{align}
\norm{u'}_{L^p(\RR)}&\leq C\abs{u'}_{L^p(\RR)}\nonumber\\
&=C \sup_{\substack{\varphi\in W^{2,p'}(\RR)\\ \varphi\not\equiv 0}}\frac{\displaystyle\int_\RR u'(x)\spt{\varphi(x)-\varphi''(x)}\D x}{\norm{\varphi}_{L^{p'}(\RR)}+\norm{\varphi''}_{L^{p'}(\RR)}}\nonumber\\
&= C\sup_{\substack{\varphi\in W^{2,p'}(\RR)\\ \varphi\not\equiv 0}}\frac{\displaystyle\int_\RR \spt{u''(x)-u(x)}\varphi'(x)\D x}{\norm{\varphi}_{L^{p'}(\RR)}+\norm{\varphi''}_{L^{p'}(\RR)}}\nonumber\\
&\leq C\norm{u''-u}_{L^p(\RR)}\sup_{\substack{\varphi\in W^{2,p'}(\RR)\\ \varphi\not\equiv 0}}\frac{\norm{\varphi'}_{L^{p'}(\RR)}}{\norm{\varphi}_{L^{p'}(\RR)}+\norm{\varphi''}_{L^{p'}(\RR)}}\label{esti1}.
\end{align}
But, from \eqref{GN-1D} for \(p\geq 2\) and \(p=q=r\), we have for \(p'>2\) (i.e., \(1<p<2\)) and all \(\varphi\in C^2_c(\RR)\)
\begin{equation}\label{esti2}
\norm{\varphi'}_{L^{p'}(\RR)}\leq C \norm{\varphi}_{L^{p'}(\RR)}^{\frac12}\norm{\varphi''}_{L^{p'}(\RR)}^{\frac12}\leq C\pt{\norm{\varphi}_{L^{p'}(\RR)}+\norm{\varphi''}_{L^{p'}(\RR)}},
\end{equation}
which by a density argument can be extended to all \(\varphi\in W^{2,p'}(\RR)\). Hence, from \eqref{esti1} and \eqref{esti2} we obtain
\[
\norm{u'}_{L^p(\RR)}\leq C\norm{u''-u}_{L^p(\RR)}\leq C\pt{\norm{u}_{L^p(\RR)}+\norm{u''}_{L^p(\RR)}}.
\]
Because this estimate is valid for all \(u\in C^\infty_c(\RR)\), we can apply it, for \(\lambda>0\), to the function \(u_\lambda(x):=u(\lambda x)\) to obtain
\[
\norm{u'}_{L^p(\RR)}\leq C\pt{\frac1\lambda\norm{u}_{L^p(\RR)}+\lambda\norm{u''}_{L^p(\RR)}},
\]
and we conclude by choosing
\[
\lambda=\norm{u}_{L^p(\RR)}^{\frac12}\norm{u''}_{L^p(\RR)}^{-\frac12}.
\]
\end{proof}

\begin{remark}
Note that when \(p\geq 2\), Proof \#1 is perhaps the simplest of the proofs given (in fact, it appears as an exercise in \cite{Eva2010}*{Section 5.10, Problem 10}). However, to prove the general case, we use tools from functional analysis, and, as we will see later in the proof of \cref{lp-norm-equiv}, also an ODE argument as in Proof \#3 below. Also, to cover the full range of \(p\)'s, one needs to handle the cases \(p=1\) and \(p=\infty\) separately.
\end{remark}

\begin{remark}\label{rem:strong-fails}
Estimate \eqref{strong-esti} fails for every \(1\leq p<2\). Indeed, consider the ``tent'' function
\[
T(x):=(1-\abs{x})_+,
\]
fix a nonnegative \(\rho\in C^\infty_c(-1,1)\) with \(\int_\RR\rho(x)\D x=1\), and for \(0<\epsilon<\frac14\) consider \(u_\epsilon:=T\ast\rho_\epsilon\in C^\infty_c(\RR)\), where \(\rho_\epsilon(x):=\epsilon^{-1}\rho(\epsilon^{-1}x)\). Since \(T'=1\) in \((-1,0)\), \(T'=-1\) in \((0,1)\), and \(\rho_\epsilon\) is supported in \((-\epsilon,\epsilon)\), we have \(\abs{u_\epsilon'}=\abs{T'\ast\rho_\epsilon}=1\) in \((-1+\epsilon,-\epsilon)\cup(\epsilon,1-\epsilon)\), hence
\[
\int_\RR\abs{u_\epsilon'(x)}^p\D x\geq 2(1-2\epsilon)\geq 1.
\]
On the other hand, \(0\leq u_\epsilon\leq 1\), and since \(T''=\delta_{-1}-2\delta_0+\delta_1\) in the sense of distributions, we have
\[
u_\epsilon''(x)=\rho_\epsilon(x+1)-2\rho_\epsilon(x)+\rho_\epsilon(x-1).
\]
This function is supported in \((-1-\epsilon,-1+\epsilon)\cup(-\epsilon,\epsilon)\cup(1-\epsilon,1+\epsilon)\), a set of measure \(6\epsilon\), and there \(\abs{u_\epsilon''}\leq 2\epsilon^{-1}\norm{\rho}_{L^\infty(\RR)}\). Therefore
\[
\int_\RR\abs{u_\epsilon(x)u_\epsilon''(x)}^{\frac p2}\D x\leq 6\epsilon\pt{\frac{2\norm{\rho}_{L^\infty(\RR)}}{\epsilon}}^{\frac p2}=6\pt{2\norm{\rho}_{L^\infty(\RR)}}^{\frac p2}\epsilon^{1-\frac p2},
\]
which tends to \(0\) as \(\epsilon\to0\) because \(p<2\).
\end{remark}

\subsection{Proof \#2 -- Another integration by parts proof}\label{sec:proof2}

In this subsection we prove \eqref{GN-1D} for \(p=q=r\) and \(1<p<2\).

\begin{proof}
We begin from the following formal integration by parts identity related to the \(L^p\) dissipativity of the Laplacian (see for instance \cite{CiMa2022}), which is valid for every \(p>1\):
\[
\int_\RR\abs{u(x)}^{p-2}u(x)u''(x)\D x
=
-\int_\RR \pt{\abs{u(x)}^{p-2}u(x)}'u'(x)\D x
=
-(p-1)\int_\RR \abs{u(x)}^{p-2}\abs{u'(x)}^2\D x,
\]
where for \(1<p<2\) the weight \(\abs{u}^{p-2}\) is singular on \(\{u=0\}\). To make this rigorous, we show that for every \(1<s\leq2\) one has
\begin{equation}\label{weighted-esti}
\int_{\{u\neq0\}}\abs{u(x)}^{s-2}\abs{u'(x)}^2\D x\leq\frac1{s-1}\int_\RR\abs{u(x)}^{s-1}\abs{u''(x)}\D x.
\end{equation}
Inequality \eqref{weighted-esti} is the case \(p=2\), \(\theta=1-\frac s2\) of \cite{KaPe2012}*{Proposition 6.1}, where it is proved for every \(s>1\); we include the short proof for completeness. For \(\epsilon>0\) let
\[
G_\epsilon(t):=\int_0^t\pt{\sigma^2+\epsilon^2}^{\frac{s-2}2}\D\sigma,\qquad t\in\RR.
\]
The function \(G_\epsilon\) is smooth and \(G_\epsilon(0)=0\), so \(G_\epsilon(u)\in C^\infty_c(\RR)\) with \(\pt{G_\epsilon(u)}'=\pt{u^2+\epsilon^2}^{\frac{s-2}2}u'\). Moreover, since \(s\leq2\) we have \(\pt{\sigma^2+\epsilon^2}^{\frac{s-2}2}\leq\sigma^{s-2}\) for \(\sigma>0\), thus
\begin{equation}\label{G-eps-bound}
\abs{G_\epsilon(t)}\leq\frac{\abs{t}^{s-1}}{s-1}\qquad\text{for all }t\in\RR.
\end{equation}
Integrating by parts we get
\begin{align*}
\int_\RR\pt{u(x)^2+\epsilon^2}^{\frac{s-2}2}\abs{u'(x)}^2\D x&=\int_\RR u'(x)\pt{G_\epsilon(u(x))}'\D x\\
&=-\int_\RR u''(x)G_\epsilon(u(x))\D x\\
&\leq\frac1{s-1}\int_\RR\abs{u(x)}^{s-1}\abs{u''(x)}\D x.
\end{align*}
Since \(s\leq 2\), the integrand on the left-hand side increases to \(\abs{u}^{s-2}\abs{u'}^2\) in \(\{u\neq0\}\) as \(\epsilon\downarrow0\). Hence, restricting the left-hand side to \(\{u\neq0\}\) and using the monotone convergence theorem, we obtain \eqref{weighted-esti}.

Now let \(1<p<2\). A point where \(u=0\) and \(u'\neq0\) is an isolated point of \(\{u=0\}\), so \(u'=0\) a.e. in \(\{u=0\}\). Therefore H{\"o}lder's inequality with exponents \(\frac2p\) and \(\frac2{2-p}\), and \eqref{weighted-esti} with \(s=p\), give
\begin{align*}
\int_\RR\abs{u'(x)}^p\D x&=\int_{\{u\neq0\}}\pt{\abs{u(x)}^{p-2}\abs{u'(x)}^2}^{\frac p2}\abs{u(x)}^{\frac{(2-p)p}2}\D x\\
&\leq\pt{\int_{\{u\neq0\}}\abs{u(x)}^{p-2}\abs{u'(x)}^2\D x}^{\frac p2}\pt{\int_\RR\abs{u(x)}^p\D x}^{\frac{2-p}2}\\
&\leq\frac1{(p-1)^{\frac p2}}\pt{\int_\RR\abs{u(x)}^{p-1}\abs{u''(x)}\D x}^{\frac p2}\pt{\int_\RR\abs{u(x)}^p\D x}^{\frac{2-p}2}.
\end{align*}
Finally, H{\"o}lder's inequality with exponents \(p'=\frac p{p-1}\) and \(p\) gives
\[
\int_\RR\abs{u(x)}^{p-1}\abs{u''(x)}\D x\leq\norm{u}_{L^p(\RR)}^{p-1}\norm{u''}_{L^p(\RR)},
\]
and therefore
\[
\norm{u'}_{L^p(\RR)}^p\leq\frac1{(p-1)^{\frac p2}}\norm{u}_{L^p(\RR)}^{\frac{p(p-1)}2+\frac{p(2-p)}2}\norm{u''}_{L^p(\RR)}^{\frac p2}=\frac1{(p-1)^{\frac p2}}\norm{u}_{L^p(\RR)}^{\frac p2}\norm{u''}_{L^p(\RR)}^{\frac p2},
\]
which is \eqref{GN-1D} for \(p=q=r\), with \(C=(p-1)^{-\frac12}\).
\end{proof}

\begin{remark}
When \(p\geq2\), Proof \#1 gives
\[
\norm{u'}_{L^p(\RR)}^2\leq(p-1)\norm{u}_{L^p(\RR)}\norm{u''}_{L^p(\RR)}.
\]
When \(1<p<2\), Proof \#2 gives the same inequality with \(p-1\) replaced by \(\frac1{p-1}=p'-1\). Hence, for every \(1<p<\infty\),
\[
\norm{u'}_{L^p(\RR)}^2\leq\pt{\max\{p,p'\}-1}\norm{u}_{L^p(\RR)}\norm{u''}_{L^p(\RR)}.
\]
\end{remark}

\begin{remark}\label{rem:general-qr}
Proofs \#1 and \#2 also give \eqref{GN-1D} directly for exponents \(p\), \(q\), \(r\) that are not all equal. When \(p\geq 2\) this follows directly from \eqref{strong-esti}: H{\"o}lder's inequality and \(\frac2p=\frac1q+\frac1r\) give
\[
\norm{uu''}_{L^{\frac p2}(\RR)}\leq\norm{u}_{L^q(\RR)}\norm{u''}_{L^r(\RR)}
\Rightarrow
\norm{u'}_{L^p(\RR)}^2\leq(p-1)\norm{u}_{L^q(\RR)}\norm{u''}_{L^r(\RR)}.
\]

Now let \(1<p<2\) and \(r>1\), let \(r'\) be such that \(\frac1r+\frac1{r'}=1\), and set \(s:=1+\frac q{r'}\). Since \(r>1\) we have \(s>1\), and since
\[
\frac1q-\frac1{r'}=\frac1q+\frac1r-1=\frac2p-1>0,
\]
we have \(s<2\). Moreover
\[
2-s=1-\frac q{r'}=q\pt{\frac1q-\frac1{r'}}=\frac{q(2-p)}p,\qquad\text{that is}\qquad\frac{(2-s)p}{2-p}=q,\qquad\text{and}\qquad(s-1)r'=q.
\]
Repeating Proof \#2 with \(\abs{u}^{s-2}\) in place of \(\abs{u}^{p-2}\), namely using H{\"o}lder's inequality with exponents \(\frac2p\) and \(\frac2{2-p}\), and \eqref{weighted-esti}, we get
\[
\int_\RR\abs{u'(x)}^p\D x\leq\frac1{(s-1)^{\frac p2}}\pt{\int_\RR\abs{u(x)}^{s-1}\abs{u''(x)}\D x}^{\frac p2}\pt{\int_\RR\abs{u(x)}^{q}\D x}^{\frac{2-p}2},
\]
while H{\"o}lder's inequality with exponents \(r'\) and \(r\) gives
\[
\int_\RR\abs{u(x)}^{s-1}\abs{u''(x)}\D x\leq\norm{u}_{L^q(\RR)}^{\frac q{r'}}\norm{u''}_{L^r(\RR)}.
\]
Since
\[
\frac p2\cdot\frac q{r'}+\frac{q(2-p)}2=\frac q2\pt{2-\frac pr}=\frac{qp}2\pt{\frac2p-\frac1r}=\frac p2,
\]
we arrive at
\[
\norm{u'}_{L^p(\RR)}^2\leq\frac{r'}q\norm{u}_{L^q(\RR)}\norm{u''}_{L^r(\RR)}.
\]
Thus Proofs \#1 and \#2 cover every admissible triple with \(p<\infty\) except those with \(r=1\) and \(1\leq p<2\) (and then \(q=\frac p{2-p}\)), which include \(p=q=r=1\). For these exponents the argument breaks down at \(s=1\): the constant \(\frac1{s-1}\) blows up, and in fact \(\int_{\{u\neq0\}}\abs{u}^{-1}\abs{u'}^2\D x=\infty\) as soon as \(u\) has a zero \(x_0\) with \(u'(x_0)\neq0\), because
\[
\abs{u'(x)}^2\abs{u(x)}^{-1}\geq\frac{\abs{u'(x_0)}}{8\abs{x-x_0}}\qquad\text{for }x\text{ close to }x_0.
\]
The case \(p=q=r=1\) is also covered by Proof \#3 below. All these exponents are covered by Nirenberg's result: he proves the inequality for \(1<r<\infty\) and obtains \(r=1\) by letting \(r\to1\) \cite{Nir59}*{p. 129}. A different proof, based on a decomposition of \(\RR\) into dyadic level sets of \(u'\), follows from \cite{LeRoSo2025}*{Theorem 1.1}, with \(X=L^1(\RR)\) and \(Y=L^q(\RR)\).
\end{remark}

\begin{prop}[\cite{Nir59}*{p. 129}, \cite{LeRoSo2025}*{Theorem 1.1}]\label{prop:r=1}
Let \(1\leq p<2\) and \(q=\frac p{2-p}\), so that \(\frac2p=\frac1q+1\). Then there exists a constant \(C>0\), depending only on \(p\), such that for every \(u\in C^\infty_c(\RR)\)
\begin{equation}\label{GN-r=1}
\norm{u'}_{L^p(\RR)}^2\leq C\norm{u}_{L^q(\RR)}\norm{u''}_{L^1(\RR)}.
\end{equation}
\end{prop}

\subsection{Proof \#3 -- ODE argument}\label{sec:proof3}

In this subsection we prove \eqref{GN-1D} for \(p=q=r\) and \(1\leq p\leq\infty\).

\begin{proof}
For \(f:=-u''+u\in C^\infty_c(\RR)\), we have that \(u\) can be considered as the unique solution of
\begin{equation}\label{ode1}
\left\{
\begin{aligned}
-u''+u&=f&&\text{in }(-\infty,\infty),\\
\lim_{\abs x\to\infty}u(x)&=0.
\end{aligned}
\right.
\end{equation}
It is not difficult to solve this equation explicitly (e.g., using variation of constants or Fourier transforms) and find that \(u\) can be written in terms of \(f\) as
\[
u(x)=\Phi\ast f(x),
\]
for
\begin{equation}\label{defn-phi}
\Phi(x):=\frac12e^{-\abs{x}}.
\end{equation}
Observe that \(\Phi\in W^{1,q}(\RR)\) for all \(1\leq q\leq\infty\) (here \(W^{1,q}\) is the usual Sobolev space), and that it satisfies \eqref{ode1} for \(f=\delta_0\) in the sense of distributions\footnote{Note that \(\hat\Phi(\xi) \sim (1+\abs{\xi}^2)^{-1}\).}. If we compute the (weak) derivative of \(u\) using this formula, we obtain
\[
u'(x)=\Phi'\ast f(x),
\]
where a direct computation shows that\footnote{Note that \(\hat{u'}\sim\xi \hat \Phi \hat f\).}
\[
\Phi'(x)=-\frac12\operatorname{sgn}(x)e^{-\abs{x}},\qquad x\neq0.
\]
In particular, \(\Phi'\in L^1(\RR)\) with \(\norm{\Phi'}_{L^1(\RR)}=1\), so Young's inequality for convolutions gives, for all \(1\leq p\leq \infty\),
\begin{equation}\label{young}
\norm{u'}_{L^p(\RR)} \leq \norm{\Phi'}_{L^1(\RR)}\norm{f}_{L^p(\RR)}=\norm{f}_{L^p(\RR)}.
\end{equation}
Since \(f=u-u''\), \eqref{young} gives
\[
\norm{u'}_{L^p(\RR)} \leq \norm{u}_{L^p(\RR)}+\norm{u''}_{L^p(\RR)}.
\]
Because the above estimate is valid for all \(u\in C^\infty_c(\RR)\), we can use it for \(\lambda>0\) and the function \(u_\lambda(x):=u(\lambda x)\) to obtain
\begin{equation}\label{ode-scaling}
\norm{u'}_{L^p(\RR)}\leq \frac{1}{\lambda}\norm{u}_{L^p(\RR)}+\lambda\norm{u''}_{L^p(\RR)},
\end{equation}
and the interpolation inequality when \(p=q=r\) follows by letting
\[
\lambda=\norm{u''}_{L^p(\RR)}^{-\frac12}\norm{u}_{L^p(\RR)}^{\frac12}.
\]
Observe that we have proven \eqref{GN-1D} for \(p=q=r\), with \(C=2\).
\end{proof}

\subsection{Proof \#4 -- \texorpdfstring{\(L^p\)}{Lp} theory}\label{sec:proof4}

In this subsection we prove \eqref{GN-1D} for \(p=q=r\) and \(1<p<\infty\).

\begin{proof}
This last proof uses the \(L^p\) theory of singular integrals, more precisely, the following classical result
\begin{prop}[\cite{Ste70}*{Proposition 3, p. 59}]\label{Lp-theory}
Let \(N\geq 2\) and \(1<p<\infty\). Then there exists a constant \(C>0\) such that for all \(v\in C^\infty_c(\RR^N)\) the following inequality holds
\begin{equation}
\norm{D^2v}_{L^p(\RR^N)}\leq C\norm{\Delta v}_{L^p(\RR^N)}.
\end{equation}
\end{prop}

To use \cref{Lp-theory} in the proof of \eqref{GN-1D} when \(p=q=r>1\), consider the function \(v(x,y):=u(x)u(y)\), which belongs to \(C^\infty_c(\RR^2)\) and satisfies
\begin{gather*}
\frac{\partial^2 v}{\partial x\partial y}(x,y)=u'(x)u'(y),\\
\Delta v(x,y)=u''(x)u(y)+u(x)u''(y).
\end{gather*}
Therefore, Fubini's theorem and \cref{Lp-theory} imply that we can write
\begin{align*}
\pt{\int_{\RR}\abs{u'(x)}^p\D x}^2 %&=\pt{\int_{\RR}\abs{u'(x)}^p\D x}\pt{\int_{\RR}\abs{u'(y)}^p\D y}\\
&=\iint_{\RR^2}\abs{u'(x)u'(y)}^p\D x\D y\\
&=\iint_{\RR^2}\abs{\frac{\partial^2 v}{\partial x\partial y}(x,y)}^p\D x\D y\\
&\leq C^p\iint_{\RR^2}\abs{\Delta v(x,y)}^p\D x\D y\\
&= C^p\iint_{\RR^2}\abs{u''(x)u(y)+u''(y)u(x)}^p\D x\D y\\
%&\leq 2^{p-1}C\iint_{\RR^2}\spt{\abs{u''(x)u(y)}^p+\abs{u''(y)u(x)}^p}\D x\D y\\
&\leq 2^{p-1}C^p\iint_{\RR^2}\abs{u''(x)u(y)}^p\D x\D y\\
&= 2^{p-1}C^p\pt{\int_{\RR}\abs{u''(x)}^p\D x}\pt{\int_{\RR}\abs{u(y)}^p\D y},
\end{align*}
which is exactly \eqref{GN-1D} for \(p=q=r\).
\end{proof}

\begin{remark}
In Proof \#4 the cases \(p=1\) and \(p=\infty\) have to be treated separately.
\end{remark}

\subsection{Proof of \texorpdfstring{\cref{gn-thm}}{Theorem 1}}\label{sec:general-a}

In this subsection we prove \cref{gn-thm} for all \(\frac12\leq a\leq1\).

\begin{proof}
The right-hand side of \eqref{gen-GN2} is an affine function of \(a\) with slope \(\frac1r-\frac1q-2<0\), and for \(a=\frac12\) it equals \(\frac12\pt{\frac1q+\frac1r}\leq1\). Hence \(\frac1p\leq1\), that is, \(p\geq1\).

\emph{Step 1: the case \(p\geq\max\{q,r\}\).} Recall from Proof \#3 that \(u'=\Phi'\ast f\) with \(f=u-u''\). A direct computation gives
\[
\norm{\Phi'}_{L^t(\RR)}=\frac12\pt{\frac2t}^{\frac1t}\leq1\quad\text{for }1\leq t<\infty,\qquad\text{and}\qquad\norm{\Phi'}_{L^\infty(\RR)}=\frac12.
\]
Since \(p\geq q\) and \(p\geq r\), there are \(t,s\in[1,\infty]\) such that \(1+\frac1p=\frac1t+\frac1q=\frac1s+\frac1r\). Hence Young's inequality for convolutions gives, as an improvement of \eqref{young},
\begin{equation}\label{esti5}
\norm{u'}_{L^p(\RR)}\leq\norm{\Phi'}_{L^t(\RR)}\norm{u}_{L^q(\RR)}+\norm{\Phi'}_{L^s(\RR)}\norm{u''}_{L^r(\RR)}\leq\norm{u}_{L^q(\RR)}+\norm{u''}_{L^r(\RR)}.
\end{equation}
Using \eqref{esti5} for the function \(u_\lambda(x):=u(\lambda x)\), \(\lambda>0\), gives
\begin{equation}\label{f-lambda}
\norm{u'}_{L^p(\RR)}\leq\lambda^{-\sigma}\norm{u}_{L^q(\RR)}+\lambda^{\tau}\norm{u''}_{L^r(\RR)},\qquad\text{where}\qquad\sigma:=1+\frac1q-\frac1p,\quad\tau:=1+\frac1p-\frac1r.
\end{equation}
Here \(\sigma\geq1\), because \(p\geq q\), and \(\tau\geq\frac1p\geq0\). Moreover, \eqref{gen-GN2} is equivalent to \(\sigma=a(\sigma+\tau)\), and hence \(\tau=(1-a)(\sigma+\tau)\). If \(\tau=0\), then \(a=1\), and letting \(\lambda\to\infty\) in \eqref{f-lambda} gives \eqref{gen-GN} with \(C=1\). If \(\tau>0\), the right-hand side of \eqref{f-lambda} attains its minimum at
\[
\underline\lambda:=\pt{\frac{\sigma\norm{u}_{L^q(\RR)}}{\tau\norm{u''}_{L^r(\RR)}}}^{\frac1{\sigma+\tau}},
\]
and a direct computation gives \eqref{gen-GN} with \(C=a^{-a}(1-a)^{a-1}\). This constant is a decreasing function of \(a\in[\frac12,1]\), so it is at most \(2\), its value for \(a=\frac12\).

\emph{Step 2: the general case.} Fix \(q\) and \(r\), and for \(\frac12\leq b\leq1\) let \(\frac1{p(b)}\) be the right-hand side of \eqref{gen-GN2} with \(a\) replaced by \(b\). This is a decreasing affine function of \(b\). For \(b=\frac12\) it is at least \(\min\{\frac1q,\frac1r\}\), and for \(b=1\) it equals \(\frac1r-1\leq0\). Hence there is \(a_1\in[\frac12,1]\) such that \(p(a_1)=\max\{q,r\}\). If \(a\geq a_1\), then \(p=p(a)\geq\max\{q,r\}\), and Step 1 applies. If \(\frac12\leq a<a_1\), we write \(a=\frac\theta2+(1-\theta)a_1\) with \(\theta\in(0,1]\). Then \(\frac1p=\frac\theta{p(\frac12)}+\frac{1-\theta}{p(a_1)}\), so H{\"o}lder's inequality gives
\[
\norm{u'}_{L^p(\RR)}\leq\norm{u'}_{L^{p(\frac12)}(\RR)}^\theta\norm{u'}_{L^{p(a_1)}(\RR)}^{1-\theta}.
\]
Since \(\frac2{p(\frac12)}=\frac1q+\frac1r\), we can estimate the first factor by \eqref{GN-1D}, and the second one by Step 1. Since \(\frac\theta2+(1-\theta)a_1=a\), we obtain \eqref{gen-GN}.
\end{proof}

\begin{remark}\label{rem:general-a}
Step 1 shows that, if \(p\geq\max\{q,r\}\), then \eqref{gen-GN} holds with the explicit constant \(C=a^{-a}(1-a)^{a-1}\leq2\), where \(0^0:=1\). The ODE argument alone gives \eqref{gen-GN} exactly in this case. For \(a=\frac12\) this means \(p=q=r\), which is the case treated in Proof \#3, with the same constant \(2\). Hence the ODE argument does not give \eqref{GN-1D} for exponents that are not all equal. On the other hand, for all \(q\) and \(r\), Step 1 directly gives the case \(p=\infty\), that is, \eqref{gen-GN} with \(\frac1p=0\) and
\[
a=\frac{1+\frac1q}{2+\frac1q-\frac1r}.
\]
When \(a=1\), \eqref{gen-GN2} forces \(r=1\) and \(p=\infty\), and \eqref{gen-GN} becomes the one-dimensional Sobolev inequality for \(u'\),
\[
\norm{u'}_{L^\infty(\RR)}\leq C\norm{u''}_{L^1(\RR)},
\]
which in fact holds with \(C=\frac12\), since \(2u'(x)=\int_{-\infty}^xu''(t)\D t-\int_x^\infty u''(t)\D t\).
\end{remark}

In higher dimensions, the endpoint estimate is Sobolev's inequality. We use the following lemma in \cref{sec:dimN,sec:weights}.

\begin{lem}\label{lem:general-a}
Let \(\mu\) be a Borel measure on \(\RR^N\), and write \(L^t:=L^t(\RR^N,\mu)\). Let \(n>1\), \(1\leq q\leq\infty\) and \(1\leq r<n\), and define \(p_0\) and \(r^*\) by
\[
\frac2{p_0}=\frac1q+\frac1r\qquad\text{and}\qquad\frac1{r^*}=\frac1r-\frac1n.
\]
Assume that \(\mu\) supports the Sobolev inequality
\begin{equation}\label{sob-lem}
\norm{v}_{L^{r^*}}\leq S\norm{\nabla v}_{L^r}\qquad\text{for every Lipschitz function }v\text{ with compact support.}
\end{equation}
Let \(u\) be a function on \(\RR^N\), and let \(g,h\geq0\) be functions on \(\RR^N\) such that \(g\) is Lipschitz with compact support, \(\abs{\nabla g}\leq h\) a.e., and
\begin{equation}\label{half-lem}
\norm{g}_{L^{p_0}}^2\leq C_0\norm{u}_{L^q}\norm{h}_{L^r}.
\end{equation}
Then, for every \(\frac12\leq a\leq1\),
\begin{equation}\label{gen-lem}
\norm{g}_{L^p}\leq C_0^{1-a}S^{2a-1}\norm{u}_{L^q}^{1-a}\norm{h}_{L^r}^a,\qquad\text{where}\qquad\frac1p=\frac1n+a\pt{\frac1r-\frac2n}+\frac{1-a}q.
\end{equation}
\end{lem}

\begin{proof}
Let \(b:=2a-1\in[0,1]\). A direct computation shows that \(\frac1p=\frac b{r^*}+\frac{1-b}{p_0}\). Hence H{\"o}lder's inequality, \eqref{sob-lem} with \(v=g\), and \eqref{half-lem} give
\[
\norm{g}_{L^p}\leq\norm{g}_{L^{r^*}}^b\norm{g}_{L^{p_0}}^{1-b}\leq\pt{S\norm{h}_{L^r}}^b\pt{C_0\norm{u}_{L^q}\norm{h}_{L^r}}^{\frac{1-b}2},
\]
which is \eqref{gen-lem}, because \(\frac{1-b}2=1-a\).
\end{proof}

\section{Functions of several variables}\label{sec:dimN}

In this section we extend \eqref{GN-1D} to functions of several variables, starting with the following result (see \cite{Nir59}*{Inequality (2.5), p. 129}).
\begin{thm}\label{gn-thm-ND}
Let \(1\leq p\leq \infty\) and \(1\leq q,r\leq\infty\) satisfy \(\frac2p=\frac1q+\frac1r\). Then there exists a constant \(C>0\) such that for every \(u\in C^\infty_c(\RR^N)\)
\begin{equation}\label{GN-ND}
\norm{\nabla u}_{L^p(\RR^N)}\leq C\norm{u}_{L^q(\RR^N)}^{\frac12}\norm{D^2u}_{L^r(\RR^N)}^{\frac12}.
\end{equation}
\end{thm}

We include the proof of \cref{gn-thm-ND} for the sake of completeness.
\begin{proof}
This result follows directly from \eqref{GN-1D} by integrating it over the ``remaining variables''. For example, if \(x=(x_1,\bar x)\in\RR\times\RR^{N-1}\) and if we want to estimate \(u_{x_1}:=\frac{\partial u}{\partial x_1}\), we apply \eqref{GN-1D} to \(u(\cdot, \bar x)\) and obtain
\[
\int_{\RR}\abs{u_{x_1}(x_1,\bar x)}^p\D x_1\leq C^p\pt{\int_{\RR}\abs{u_{x_1x_1}(x_1,\bar x)}^r\D x_1}^{\frac p{2r}}\pt{\int_{\RR}\abs{u(x_1,\bar x)}^q\D x_1}^{\frac p{2q}}.
\]
Since this holds for all \(\bar x\), we can integrate over \(\bar x\in \RR^{N-1}\) and use H{\"o}lder's inequality with \(1=\frac1s+\frac1t\) to obtain
\begin{multline*}
\int_{\RR^{N-1}}\int_{\RR}\abs{u_{x_1}(x_1,\bar x)}^p\D x_1\D \bar x\leq C^p\int_{\RR^{N-1}}\pt{\int_{\RR}\abs{u_{x_1x_1}(x_1,\bar x)}^r\D x_1}^{\frac p{2r}}\pt{\int_{\RR}\abs{u(x_1,\bar x)}^q\D x_1}^{\frac p{2q}}\D \bar x\\
\leq C^p\pt{\int_{\RR^{N-1}}\pt{\int_{\RR}\abs{u_{x_1x_1}(x_1,\bar x)}^r\D x_1}^{\frac p{2r}\cdot s}\D \bar x}^{\frac 1s}\pt{\int_{\RR^{N-1}}\pt{\int_{\RR}\abs{u(x_1,\bar x)}^q\D x_1}^{\frac p{2q}\cdot t}\D \bar x}^{\frac1t},
\end{multline*}
and choosing \(s=\frac{2r}p\) and \(t=\frac{2q}p\) (observe that \(\frac2p=\frac1q+\frac1r\) implies that \(s,t\geq 1\)), we obtain
\[
\int_{\RR^N}\abs{u_{x_1}(x)}^p\D x\leq C^p\pt{\int_{\RR^N}\abs{u_{x_1x_1}(x)}^r\D x}^{\frac p{2r}}\pt{\int_{\RR^N}\abs{u(x)}^q\D x}^{\frac p{2q}}.
\]
If \(p<\infty\) and \(q=\infty\) (respectively \(r=\infty\)), then \(\frac p{2r}=1\) (respectively \(\frac p{2q}=1\)), and instead of H{\"o}lder's inequality it suffices to bound \(\norm{u(\cdot,\bar x)}_{L^\infty(\RR)}\leq\norm{u}_{L^\infty(\RR^N)}\) (respectively \(\norm{u_{x_1x_1}(\cdot,\bar x)}_{L^\infty(\RR)}\leq\norm{u_{x_1x_1}}_{L^\infty(\RR^N)}\)) before integrating in \(\bar x\). If \(p=\infty\), then \(q=r=\infty\), and it suffices to take the supremum over \(\bar x\) in \eqref{GN-1D}.
Since the above can be repeated for each \(x_i\), \(i=2,\ldots,N\), we obtain \eqref{GN-ND}.
\end{proof}

\begin{remark}\label{rem:general-ND}
The slicing argument above does not extend to the general form \eqref{gen-GN} of \cref{gn-thm}. On each line it requires the one-dimensional relation \eqref{gen-GN2}, while H{\"o}lder's inequality in the remaining variables requires \(\frac1p=\frac ar+\frac{1-a}q\), and both hold only when \(a=\frac12\). In \(\RR^N\), the dilations \(u\mapsto u(\lambda\,\cdot)\) lead instead to the relation
\begin{equation}\label{gen-GN2-N}
\frac1p=\frac1N+a\pt{\frac1r-\frac2N}+\frac{1-a}q,
\end{equation}
which is the one in Nirenberg's theorem \cite{Nir59}*{p. 125}. If \(N\geq2\) and \(1\leq r<N\), the corresponding inequality
\begin{equation}\label{gen-GN-N}
\norm{\nabla u}_{L^p(\RR^N)}\leq C\norm{u}_{L^q(\RR^N)}^{1-a}\norm{D^2u}_{L^r(\RR^N)}^a,\qquad1\leq q\leq\infty,\quad\frac12\leq a\leq1,
\end{equation}
follows from \cref{lem:general-a}, with \(\mu\) the Lebesgue measure, \(n=N\), \(g=\abs{\nabla u}\) and \(h=\abs{D^2u}\). Indeed, \eqref{sob-lem} is Sobolev's inequality, Kato's inequality gives \(\abs{\nabla\abs{\nabla u}}\leq\abs{D^2u}\) a.e., and \eqref{half-lem} is \cref{gn-thm-ND} with exponents \(p_0\), \(q\) and \(r\), where \(p_0\) is as in \cref{lem:general-a}.
\end{remark}

\begin{remark}\label{rem:general-kernel}
Proof \#3, the ODE argument, also extends to \(\RR^N\), with the Bessel kernel in place of the fundamental solution \(\Phi\) of \(-u''+u=f\). In this form it gives the general form \eqref{gen-GN-N} in some of the cases \(r>N\), which are not covered by \cref{rem:general-ND}. Let \(G_1\) be the kernel of \((1-\Delta)^{-1}\) in \(\RR^N\), so that \(G_1=\Phi\) if \(N=1\), and for \(\lambda>0\) let \(G_\lambda(x):=\lambda^{N-2}G_1(\lambda x)\), the kernel of \((\lambda^2-\Delta)^{-1}\). Then \(\nabla u=\nabla G_\lambda\ast(\lambda^2u-\Delta u)\), and \(\nabla G_1\in L^t(\RR^N)\) for \(1\leq t<\frac N{N-1}\), with
\[
\norm{\nabla G_\lambda}_{L^t(\RR^N)}=\lambda^{N-1-\frac Nt}\norm{\nabla G_1}_{L^t(\RR^N)}.
\]
If \(q>N\) and \(r>N\), Young's inequality with exponents \(q'\) and \(r'\), where \(\frac1q+\frac1{q'}=\frac1r+\frac1{r'}=1\), gives
\[
\norm{\nabla u}_{L^\infty(\RR^N)}\leq C\pt{\lambda^{1+\frac Nq}\norm{u}_{L^q(\RR^N)}+\lambda^{\frac Nr-1}\norm{\Delta u}_{L^r(\RR^N)}}.
\]
Optimizing in \(\lambda\), as in Step 1 of the proof of \cref{gn-thm}, we obtain \eqref{gen-GN-N} with \(\frac1p=0\), that is, with
\[
a=\frac{1+\frac Nq}{2+\frac Nq-\frac Nr},
\]
and with \(\Delta u\) in place of \(D^2u\). Since \(\abs{\Delta u}\leq\sqrt N\abs{D^2u}\), H{\"o}lder's inequality, as in Step 2 of that proof, interpolates between this case and \cref{gn-thm-ND}, and gives \eqref{gen-GN-N} for all \(\frac12\leq a\leq1\) with \(\frac1p\geq0\). The remaining cases, \(r=N\), or \(r>N\) and \(q\leq N\), are not covered by these arguments. They are covered by Nirenberg's theorem \cite{Nir59}*{p. 125}, except for \(a=1\) and \(r=N\geq2\), which is his exceptional case: there the inequality fails.
\end{remark}

Proofs \#1 and \#2 can also be carried out directly in \(\RR^N\), and they give the following theorem, which is \cref{thmA} of the introduction.

\begin{thm}\label{gn-thm-RN}
Let \(N\geq1\), and let \(1<p<\infty\) and \(1\leq q,r\leq\infty\) satisfy \(\frac2p=\frac1q+\frac1r\). Then for every \(u\in C^\infty_c(\RR^N)\) the following holds:
\begin{enumerate}
\item[(a)] if \(p\geq2\), then
\begin{equation}\label{GN-RN-hess}
\norm{\nabla u}_{L^p(\RR^N)}^2\leq\sqrt{N-1+(p-1)^2}\,\norm{u}_{L^q(\RR^N)}\norm{D^2u}_{L^r(\RR^N)};
\end{equation}
\item[(b)] if \(p\leq2\) and \(r>1\), then
\begin{equation}\label{GN-RN-lap}
\norm{\nabla u}_{L^p(\RR^N)}^2\leq\frac{r'}{q}\norm{u}_{L^q(\RR^N)}\norm{\Delta u}_{L^r(\RR^N)},
\end{equation}
where \(\frac1r+\frac1{r'}=1\).
\end{enumerate}
\end{thm}

\begin{proof}
If \(p=2\), then \(q=r'\), and both inequalities follow from an integration by parts, H{\"o}lder's inequality and the pointwise estimate \(\abs{\Delta u}\leq\sqrt N\abs{D^2u}\), since
\[
\int_{\RR^N}\abs{\nabla u}^2\D x=-\int_{\RR^N}u\,\Delta u\D x\leq\norm{u}_{L^q(\RR^N)}\norm{\Delta u}_{L^r(\RR^N)}.
\]
So we assume \(p\neq2\).

\emph{Part (a) is Proof \#1 in \(\RR^N\).} Since \(p>2\), the map \(\xi\mapsto\abs\xi^{p-2}\xi\) is of class \(C^1\) on \(\RR^N\), and its differential vanishes at \(\xi=0\). Hence \(X:=\abs{\nabla u}^{p-2}\nabla u\) is a compactly supported \(C^1\) vector field, and \(\operatorname{div}X=0\) on \(\{\nabla u=0\}\). In the open set \(\{\nabla u\neq0\}\), let \(\nu:=\frac{\nabla u}{\abs{\nabla u}}\), and at each point complete \(e_1:=\nu\) to an orthonormal basis \(\{e_1,\ldots,e_N\}\) of \(\RR^N\). Then
\[
\operatorname{div}X=\abs{\nabla u}^{p-2}\pt{\Delta u+(p-2)D^2u(\nu,\nu)}=\abs{\nabla u}^{p-2}\pt{(p-1)D^2u(e_1,e_1)+\sum_{i=2}^ND^2u(e_i,e_i)},
\]
and the Cauchy--Schwarz inequality gives the pointwise estimate
\[
\abs{\operatorname{div}X}\leq c_{N,p}\abs{\nabla u}^{p-2}\abs{D^2u}\quad\text{in }\RR^N,\qquad\text{where }c_{N,p}:=\sqrt{N-1+(p-1)^2}.
\]
Since \(\operatorname{div}(uX)=u\operatorname{div}X+\abs{\nabla u}^p\), the divergence theorem and H{\"o}lder's inequality with exponents \(\frac p{p-2}\) and \(\frac p2\) give
\begin{align*}
\int_{\RR^N}\abs{\nabla u}^p\D x=-\int_{\RR^N}u\operatorname{div}X\D x&\leq c_{N,p}\int_{\RR^N}\abs{\nabla u}^{p-2}\abs{u\,D^2u}\D x\\
&\leq c_{N,p}\pt{\int_{\RR^N}\abs{\nabla u}^p\D x}^{1-\frac2p}\pt{\int_{\RR^N}\abs{u\,D^2u}^{\frac p2}\D x}^{\frac2p},
\end{align*}
which is the analogue of \eqref{strong-esti} in \(\RR^N\). Finally, H{\"o}lder's inequality and \(\frac2p=\frac1q+\frac1r\) give
\[
\norm{u\,D^2u}_{L^{\frac p2}(\RR^N)}\leq\norm{u}_{L^q(\RR^N)}\norm{D^2u}_{L^r(\RR^N)},
\]
and \eqref{GN-RN-hess} follows.

\emph{Part (b) is Proof \#2 in \(\RR^N\).} Let \(1<s\leq2\), and let \(G_\epsilon\) be as in Proof \#2. Then \(G_\epsilon(u)\in C^\infty_c(\RR^N)\) and \(\nabla\pt{G_\epsilon(u)}=\pt{u^2+\epsilon^2}^{\frac{s-2}2}\nabla u\), so an integration by parts and \eqref{G-eps-bound} give
\[
\int_{\RR^N}\pt{u^2+\epsilon^2}^{\frac{s-2}2}\abs{\nabla u}^2\D x=-\int_{\RR^N}G_\epsilon(u)\,\Delta u\D x\leq\frac1{s-1}\int_{\RR^N}\abs{u}^{s-1}\abs{\Delta u}\D x,
\]
and, exactly as in Proof \#2, the monotone convergence theorem yields
\begin{equation}\label{weighted-esti-RN}
\int_{\{u\neq0\}}\abs{u}^{s-2}\abs{\nabla u}^2\D x\leq\frac1{s-1}\int_{\RR^N}\abs{u}^{s-1}\abs{\Delta u}\D x.
\end{equation}
Moreover, by the implicit function theorem, \(\{u=0\}\cap\{\nabla u\neq0\}\) is a smooth hypersurface, so it has measure zero, that is, \(\nabla u=0\) a.e.\ in \(\{u=0\}\). Proof \#2 and \cref{rem:general-qr} use only these two facts; the rest of the argument consists of two applications of H{\"o}lder's inequality. Hence, repeating that argument with \(s=1+\frac q{r'}\), and with \eqref{weighted-esti-RN} in place of \eqref{weighted-esti}, we obtain \eqref{GN-RN-lap}.
\end{proof}

\begin{remark}\label{rem:RN-consts}
For \(N=1\) the constants in \eqref{GN-RN-hess} and \eqref{GN-RN-lap} are those given by Proofs \#1 and \#2. For \(p\geq2\), G{\"u}neysu and Pigola \cite{GuPi2019}*{Theorem 1} obtain, from the same integration by parts,
\[
\norm{\nabla u}_{L^p(\RR^N)}^2\leq\norm{u}_{L^q(\RR^N)}\pt{\norm{\Delta u}_{L^r(\RR^N)}+(p-2)\norm{D^2u}_{L^r(\RR^N)}},
\]
and hence the constant \(\sqrt N+p-2\) in \eqref{GN-RN-hess}. The Cauchy--Schwarz step in the proof of part (a) improves it to \(c_{N,p}=\sqrt{N-1+(p-1)^2}\), since \((\sqrt N+p-2)^2-c_{N,p}^2=2(p-2)(\sqrt N-1)\geq0\). We have not found the constant \(\frac{r'}q\) of \eqref{GN-RN-lap} in the literature: Coulhon and Duong \cite{CoDu2003}*{Theorems 4.1 and 4.3} obtain this inequality for \(p\leq2\) with implicit constants that do not depend on the dimension, and a related estimate is stated in the recent preprint \cite{Ouh2024}. \Cref{gn-thm-RN} gives a proof of \cref{gn-thm-ND} that does not rely on the one-dimensional result. It covers every admissible triple with \(1<p<\infty\), except those with \(1<p<2\) and \(r=1\), and for \(p\leq2\) it gives the stronger inequality with \(\Delta u\) in place of \(D^2u\).

The case \(p=q=r=\infty\) follows from Landau's inequality \eqref{landau}: if \(\nabla u(x)\neq0\) and \(\nu:=\nabla u(x)/\abs{\nabla u(x)}\), then \(v(t):=u(x+t\nu)\) satisfies \(v'(0)=\abs{\nabla u(x)}\) and \(\abs{v''}\leq\norm{D^2u}_{L^\infty(\RR^N)}\), so that \(\norm{\nabla u}_{L^\infty(\RR^N)}^2\leq2\norm{u}_{L^\infty(\RR^N)}\norm{D^2u}_{L^\infty(\RR^N)}\).
\end{remark}

\begin{remark}\label{rem:RN-r1}
In contrast with the one-dimensional case (see \cref{prop:r=1}), the assumption \(r>1\) in \cref{gn-thm-RN}(b) cannot be removed when \(N\geq3\). Let \(p=\frac N{N-1}\), \(q=\frac N{N-2}\) and \(r=1\), so that \(\frac2p=\frac1q+\frac1r\). For \(0<\epsilon<\frac12\), the function \(w_\epsilon(x):=\pt{\abs x^2+\epsilon^2}^{\frac{2-N}2}\) satisfies
\[
\abs{\nabla w_\epsilon(x)}=(N-2)\abs x\pt{\abs x^2+\epsilon^2}^{-\frac N2},
\]
and
\[
\Delta w_\epsilon(x)=-N(N-2)\epsilon^2\pt{\abs x^2+\epsilon^2}^{-\frac{N+2}2},
\]
hence
\[
\norm{\Delta w_\epsilon}_{L^1(\RR^N)}=N(N-2)\int_{\RR^N}\pt{\abs y^2+1}^{-\frac{N+2}2}\D y
\]
does not depend on \(\epsilon\). Fix \(\chi\in C^\infty_c(\RR^N)\) with \(0\leq\chi\leq1\), \(\chi=1\) in \(B_1\) and \(\chi=0\) outside \(B_2\), and let \(u_\epsilon:=\chi w_\epsilon\). Since \(w_\epsilon\) and \(\nabla w_\epsilon\) are bounded in \(B_2\setminus B_1\) uniformly in \(\epsilon\), we have \(\norm{\Delta u_\epsilon}_{L^1(\RR^N)}\leq C\). On the other hand, \(\frac{(N-2)q}2=\frac N2\), \((N-1)p=N\), so a direct computation in polar coordinates gives
\[
\norm{u_\epsilon}_{L^q(\RR^N)}^q\leq\int_{B_2}\frac{\D x}{\pt{\abs x^2+\epsilon^2}^{\frac N2}}\leq C\log\frac1\epsilon\qquad\text{and}\qquad\norm{\nabla u_\epsilon}_{L^p(\RR^N)}^p\geq\int_{B_1\setminus B_\epsilon}\abs{\nabla w_\epsilon(x)}^p\D x\geq c\log\frac1\epsilon.
\]
Since \(\frac2p-\frac1q=1\), we conclude that
\[
\frac{\norm{\nabla u_\epsilon}_{L^p(\RR^N)}^2}{\norm{u_\epsilon}_{L^q(\RR^N)}\norm{\Delta u_\epsilon}_{L^1(\RR^N)}}\geq c\log\frac1\epsilon\longrightarrow\infty\qquad\text{as }\epsilon\to0.
\]
This agrees with \cref{gn-thm-ND}, since \(\norm{D^2u_\epsilon}_{L^1(\RR^N)}\) is of order \(\log\frac1\epsilon\). Observe the similarity with \cref{rem:strong-fails}.
\end{remark}

\section{Weighted measures}\label{sec:weights}

In this section we replace the Lebesgue measure on \(\RR^N\) by a measure \(w(x)\D x\), where \(w\) is a weight, that is, a locally integrable function that is positive a.e. Our main motivation is homogeneous weights, such as the power weights \(\abs x^\alpha\) of the Caffarelli--Kohn--Nirenberg inequalities \cite{CKN1984} and the monomial weights \(\abs{x_1}^{A_1}\cdots\abs{x_N}^{A_N}\) of Cabr{\'e} and Ros-Oton \cite{CR-O2013-2}. Writing \(L^p(w):=L^p(\RR^N,w\D x)\), we look for inequalities of the form
\begin{equation}\label{GN-w}
\norm{\nabla u}_{L^p(w)}^2\leq C\norm{u}_{L^q(w)}\norm{\mathcal Du}_{L^r(w)},\qquad\frac2p=\frac1q+\frac1r,
\end{equation}
where \(\mathcal Du\) is either \(\Delta u\) or the Hessian \(D^2u\). If \(w\) is homogeneous of degree \(\alpha\), that is, \(w(tx)=t^\alpha w(x)\) for \(t>0\), then \eqref{GN-w} is invariant under the dilations \(u\mapsto u(\lambda\,\cdot)\), exactly as in the unweighted case.

Integrating by parts with respect to \(w\D x\) produces the operator
\[
L_wu:=w^{-1}\operatorname{div}(w\nabla u)=\Delta u+\nabla\log w\cdot\nabla u
\]
instead of \(\Delta u\). With \(L_w\) in place of \(\Delta\), the proof of \cref{gn-thm-RN} works for every smooth positive weight. To obtain \eqref{GN-w} with the Euclidean operators, the first-order term \(\nabla\log w\cdot\nabla u\) has to be controlled.

\begin{thm}\label{gn-thm-hom}
Let \(w\) be positive and of class \(C^1\) in \(\RR^N\setminus\{0\}\), and homogeneous of degree \(\alpha\in\RR\). Let \(N_\alpha:=N+\alpha\) and \(\kappa_w:=\max_{\abs x=1}\abs{\nabla\log w(x)}\). Let \(1<p<\infty\) and \(1\leq q,r\leq\infty\) satisfy \(\frac2p=\frac1q+\frac1r\), and assume \(r<N_\alpha\). Then for every \(u\in C^\infty_c(\RR^N)\) the following holds:
\begin{enumerate}
\item[(a)] if \(p\geq2\), then
\[
\norm{\nabla u}_{L^p(w)}^2\leq\pt{\sqrt{N-1+(p-1)^2}+\frac{\kappa_wr}{N_\alpha-r}}\norm{u}_{L^q(w)}\norm{D^2u}_{L^r(w)};
\]
\item[(b)] if \(p\leq2\) and \(r>1\), then
\[
\norm{\nabla u}_{L^p(w)}^2\leq\frac{r'}q\pt{\sqrt N+\frac{\kappa_wr}{N_\alpha-r}}\norm{u}_{L^q(w)}\norm{D^2u}_{L^r(w)}.
\]
\end{enumerate}
The same holds, with \(\abs{N_\alpha-r}\) in place of \(N_\alpha-r\), for every \(u\in C^\infty_c(\RR^N\setminus\{0\})\), every \(\alpha\in\RR\) and every \(r\in[1,\infty)\setminus\{N_\alpha\}\).
\end{thm}

\begin{proof}
We first prove a Hardy inequality along rays: if \(1\leq r<N_\alpha\), then every Lipschitz function \(f\) with compact support in \(\RR^N\) satisfies
\begin{equation}\label{hardy-rays}
\int_{\RR^N}\frac{\abs f^r}{\abs x^r}\,w\D x\leq\pt{\frac r{N_\alpha-r}}^r\int_{\RR^N}\abs{\partial_\rho f}^r\,w\D x,\qquad\text{where }\partial_\rho f(x):=\frac x{\abs x}\cdot\nabla f(x).
\end{equation}
In polar coordinates \(x=\rho\theta\), the homogeneity of \(w\) gives \(w(x)\D x=\rho^{N_\alpha-1}w(\theta)\D\rho\D\sigma(\theta)\), so it suffices to prove, for \(g(\rho):=f(\rho\theta)\),
\[
\int_0^\infty\abs g^r\rho^{N_\alpha-1-r}\D\rho\leq\pt{\frac r{N_\alpha-r}}^r\int_0^\infty\abs{g'}^r\rho^{N_\alpha-1}\D\rho.
\]
Since \(N_\alpha>r\), an integration by parts, with no boundary terms, and H{\"o}lder's inequality give
\begin{align*}
\int_0^\infty\abs g^r\rho^{N_\alpha-1-r}\D\rho&=-\frac r{N_\alpha-r}\int_0^\infty\abs g^{r-2}g\,g'\rho^{N_\alpha-r}\D\rho\\
&\leq\frac r{N_\alpha-r}\pt{\int_0^\infty\abs g^r\rho^{N_\alpha-1-r}\D\rho}^{1-\frac1r}\pt{\int_0^\infty\abs{g'}^r\rho^{N_\alpha-1}\D\rho}^{\frac1r},
\end{align*}
which proves the claim (for \(r=1\), \(\abs g^{r-2}g\) is the sign of \(g\)). If \(f\) vanishes near the origin, the same computation works for every \(r\neq N_\alpha\), with \(\abs{N_\alpha-r}\) in place of \(N_\alpha-r\). We apply \eqref{hardy-rays} to the Lipschitz function \(f=\abs{\nabla u}\), for which \(\abs{\partial_\rho f}\leq\abs{D^2u}\) a.e., and obtain
\begin{equation}\label{hardy-grad}
\norm{\abs x^{-1}\nabla u}_{L^r(w)}\leq\frac r{N_\alpha-r}\norm{D^2u}_{L^r(w)}.
\end{equation}

Next, \(\nabla w\) is homogeneous of degree \(\alpha-1\), and \(\abs{\nabla w(x)}\leq\kappa_w\abs x^{-1}w(x)\). Since \(N_\alpha>r\geq1\), both \(w\) and \(\nabla w\) are locally integrable. Moreover, \(w\in W^{1,1}_{\mathrm{loc}}(\RR^N)\), because when one integrates by parts against a test function in \(\RR^N\setminus B_\delta\), the boundary term on \(\partial B_\delta\) is \(O(\delta^{N_\alpha-1})\). Hence we may integrate by parts against \(w\) as in Proofs \#1 and \#2. In Proof \#1, let \(X=\abs{\nabla u}^{p-2}\nabla u\) as in the proof of \cref{gn-thm-RN}(a). Since \(\operatorname{div}(wX)=w\operatorname{div}X+\nabla w\cdot X\), we get
\[
\int_{\RR^N}\abs{\nabla u}^pw\D x=-\int_{\RR^N}u\operatorname{div}(wX)\D x\leq\int_{\RR^N}\abs u\abs{\nabla u}^{p-2}\pt{c_{N,p}\abs{D^2u}+\kappa_w\frac{\abs{\nabla u}}{\abs x}}w\D x,
\]
and H{\"o}lder's inequality as in Proof \#1 gives
\[
\norm{\nabla u}_{L^p(w)}^2\leq\norm{u}_{L^q(w)}\pt{c_{N,p}\norm{D^2u}_{L^r(w)}+\kappa_w\norm{\abs x^{-1}\nabla u}_{L^r(w)}}.
\]
In Proof \#2, since \(\operatorname{div}(w\nabla u)=w\pt{\Delta u+\nabla\log w\cdot\nabla u}\), the argument gives \eqref{weighted-esti} for the measure \(w\D x\), with \(\abs{u''}\) replaced by \(F:=\abs{\Delta u}+\kappa_w\abs x^{-1}\abs{\nabla u}\). The rest of the argument, as in \cref{rem:general-qr}, then gives
\[
\norm{\nabla u}_{L^p(w)}^2\leq\frac{r'}q\norm{u}_{L^q(w)}\norm{F}_{L^r(w)}.
\]
In both cases, \eqref{hardy-grad} and \(\abs{\Delta u}\leq\sqrt N\abs{D^2u}\) conclude the proof. If \(u\in C^\infty_c(\RR^N\setminus\{0\})\), no integrability of \(w\) near the origin is needed, and \eqref{hardy-grad} holds for every \(r\neq N_\alpha\), with \(\abs{N_\alpha-r}\) in place of \(N_\alpha-r\).
\end{proof}

\begin{remark}
Observe that if \(w=\abs x^\alpha\), then
\[
\kappa_w=\abs\alpha\max_{\abs x=1}\abs{\nabla\log\abs x}=\abs\alpha.
\]
On the other hand, \cref{gn-thm-hom} does not apply to the monomial weights \(w(x)=\abs{x_1}^{A_1}\cdots\abs{x_N}^{A_N}\) with \(A\neq0\). Indeed, they are not positive and of class \(C^1\) in \(\RR^N\setminus\{0\}\), and \(\kappa_w=\infty\), since \(\nabla\log w(x)=\sum_iA_ie_i/x_i\).
\end{remark}

\begin{thm}\label{gn-thm-superh}
Let \(w\in W^{1,1}_{\mathrm{loc}}(\RR^N)\) be positive a.e. and superharmonic, that is, \(\Delta w\leq0\) in the sense of distributions. Let \(1<p\leq2\), \(1\leq q\leq\infty\) and \(1<r\leq\infty\) satisfy \(\frac2p=\frac1q+\frac1r\). Then for every \(u\in C^\infty_c(\RR^N)\)
\[
\norm{\nabla u}_{L^p(w)}^2\leq\frac{r'}q\norm{u}_{L^q(w)}\norm{\Delta u}_{L^r(w)}.
\]
\end{thm}

\begin{proof}
For \(p=2\) we have
\[
\int_{\RR^N}\abs{\nabla u}^2w\D x=-\int_{\RR^N}u\,\Delta u\,w\D x-\frac12\int_{\RR^N}\nabla\pt{u^2}\cdot\nabla w\D x,
\]
and the last term equals \(\frac12\langle\Delta w,u^2\rangle\leq0\). For \(p<2\) we follow Proof \#2 and \cref{rem:general-qr}. With \(s:=1+\frac q{r'}\) and \(G_\epsilon\) as there, let \(H_\epsilon(t):=\int_0^tG_\epsilon(\tau)\D\tau\). Then \(G_\epsilon(u)\nabla u=\nabla\pt{H_\epsilon(u)}\), and \(H_\epsilon(u)\in C^\infty_c(\RR^N)\) is nonnegative, because \(G_\epsilon\) is odd and increasing. Therefore
\begin{align*}
\int_{\RR^N}\pt{u^2+\epsilon^2}^{\frac{s-2}2}\abs{\nabla u}^2w\D x&=-\int_{\RR^N}G_\epsilon(u)\,\Delta u\,w\D x-\int_{\RR^N}\nabla w\cdot\nabla\pt{H_\epsilon(u)}\D x\\
&\leq\frac1{s-1}\int_{\RR^N}\abs{u}^{s-1}\abs{\Delta u}w\D x,
\end{align*}
since
\[
-\int_{\RR^N}\nabla w\cdot\nabla\pt{H_\epsilon(u)}\D x=\langle\Delta w,H_\epsilon(u)\rangle\leq0.
\]
This gives \eqref{weighted-esti} for the measure \(w\D x\), with \(u''\) replaced by \(\Delta u\), and the rest of Proof \#2 and of \cref{rem:general-qr} applies without changes.
\end{proof}

\begin{remark}\label{rem:superh}
Nonconstant positive superharmonic functions on \(\RR^N\) exist only when \(N\geq3\). For instance, \(\abs x^\alpha\) is superharmonic in \(\RR^N\) if and only if \(2-N\leq \alpha\leq0\), since \(\Delta\abs x^\alpha=\alpha(N+\alpha-2)\abs x^{\alpha-2}\) for \(x\neq0\), and \(\Delta\abs x^{2-N}\) is a negative multiple of the Dirac delta. Other examples are \(\min\{1,\abs x^{2-N}\}\) and, more generally, Newtonian potentials of positive measures, which need not be homogeneous. If \(w\) is subharmonic, for example \(w=\abs x^\alpha\) with \(\alpha>0\), then the term \(\langle\Delta w,H_\epsilon(u)\rangle\) has the wrong sign, and \cref{gn-thm-hom} can be used instead.
\end{remark}

Let
\[
Mf(x):=\sup_{\rho>0}\frac1{\abs{B_\rho(x)}}\int_{B_\rho(x)}\abs{f(y)}\D y
\]
be the Hardy--Littlewood maximal function. For \(1<s<\infty\), a weight \(w\) belongs to the Muckenhoupt class \(A_s\) if
\[
[w]_{A_s}:=\sup_B\pt{\frac1{\abs B}\int_Bw\D x}\pt{\frac1{\abs B}\int_Bw^{-\frac1{s-1}}\D x}^{s-1}<\infty,
\]
where the supremum is taken over all balls \(B\subset\RR^N\). It belongs to \(A_1\) if \(Mw\leq[w]_{A_1}w\) a.e. for some constant \([w]_{A_1}\). By Muckenhoupt's theorem \cite{Muck1972}, if \(1<s<\infty\), then \(M\) is bounded on \(L^s(w)\) if and only if \(w\in A_s\). We denote its norm by \(\norm{M}_{L^s(w)}\), and we set \(\norm{M}_{L^\infty(w)}:=1\). Recall that \(A_s\subset A_t\) if \(s\leq t\), that \(\abs x^\alpha\in A_s\) if and only if \(-N<\alpha<N(s-1)\), and that \(\abs x^\alpha\in A_1\) if and only if \(-N<\alpha\leq0\).

\begin{lem}\label{lem:pointwise}
Let
\[
\beta_N:=\frac{\sqrt\pi\,\Gamma\pt{\frac{N+1}2}}{\Gamma\pt{\frac N2}},
\]
so that \(\beta_1=1\). Then for every \(u\in C^\infty_c(\RR^N)\) and every \(x\in\RR^N\)
\begin{equation}\label{pointwise}
\abs{\nabla u(x)}^2\leq4\beta_N^2\,Mu(x)\,M(\Delta u)(x).
\end{equation}
\end{lem}

\begin{proof}
This is Proof \#3 in pointwise form. Let \(G_1\) and \(G_\lambda\) be as in \cref{rem:general-kernel}. As in Proof \#3, \(u=G_\lambda\ast f\) with \(f:=\lambda^2u-\Delta u\), hence \(\nabla u=\nabla G_\lambda\ast f\). The function \(\abs{\nabla G_1}\) is radial, decreasing and integrable, with \(\norm{\nabla G_1}_{L^1(\RR^N)}=\beta_N\).\footnote{In terms of the modified Bessel functions \(K_\nu\),
\[
G_1(x)=(2\pi)^{-\frac N2}\abs x^{1-\frac N2}K_{\frac N2-1}(\abs x)\qquad\text{and}\qquad\abs{\nabla G_1(x)}=(2\pi)^{-\frac N2}\abs x^{1-\frac N2}K_{\frac N2}(\abs x).
\]
The latter is decreasing in \(\abs x\), because
\[
\frac{\D}{\D t}\pt{t^{1-\frac N2}K_{\frac N2}(t)}=-t^{-\frac N2}\pt{(N-1)K_{\frac N2}(t)+tK_{\frac N2-1}(t)},
\]
and \(\beta_N\) follows from
\[
\int_0^\infty t^{\frac N2}K_{\frac N2}(t)\D t=2^{\frac N2-1}\Gamma\pt{\tfrac12}\Gamma\pt{\tfrac{N+1}2}.
\]
For instance, if \(N=3\), then \(\abs{\nabla G_1(x)}=(1+\abs x)e^{-\abs x}/(4\pi\abs x^2)\) and \(\beta_3=2\).} Hence \(\abs{\nabla G_\lambda}\) is radial and decreasing, and \(\norm{\nabla G_\lambda}_{L^1(\RR^N)}=\frac{\beta_N}\lambda\). For such kernels we have \cite{Ste70}*{Chapter III, \S2.2}
\[
\abs{\nabla G_\lambda}\ast\abs f\leq\norm{\nabla G_\lambda}_{L^1(\RR^N)}Mf.
\]
Therefore
\[
\abs{\nabla u(x)}\leq\frac{\beta_N}\lambda M\pt{\lambda^2u-\Delta u}(x)\leq\beta_N\pt{\lambda\,Mu(x)+\frac1\lambda M(\Delta u)(x)}.
\]
This is a pointwise version of inequality \eqref{ode-scaling} of Proof \#3, with \(\lambda\) replaced by \(\frac1\lambda\). Since it holds at each point \(x\), the parameter \(\lambda\) may now depend on \(x\). Since \(u\not\equiv0\), we have \(Mu>0\) and \(M(\Delta u)>0\) everywhere, and choosing \(\lambda=\pt{M(\Delta u)(x)/Mu(x)}^{\frac12}\) we obtain \eqref{pointwise}.
\end{proof}

The following theorem contains \cref{thmB} of the introduction.

\begin{thm}\label{gn-thm-Ap}
Let \(\beta_N\) be as in \cref{lem:pointwise}, let \(1<q,r\leq\infty\), let \(1<p\leq\infty\) be given by \(\frac2p=\frac1q+\frac1r\), and let \(w\in A_{\min\{q,r\}}\) (no condition on \(w\) is needed if \(q=r=\infty\)). Then for every \(u\in C^\infty_c(\RR^N)\)
\begin{equation}\label{GN-Ap}
\norm{\nabla u}_{L^p(w)}^2\leq4\beta_N^2\norm{M}_{L^q(w)}\norm{M}_{L^r(w)}\norm{u}_{L^q(w)}\norm{\Delta u}_{L^r(w)}.
\end{equation}
Moreover, if \(w\in A_1\), then for every \(u\in C^\infty_c(\RR^N)\)
\[
\norm{\nabla u}_{L^1(w)}^2\leq4\beta_N^2[w]_{A_1}^2\norm{u}_{L^1(w)}\norm{\Delta u}_{L^1(w)}.
\]
\end{thm}

\begin{proof}
If \(p=\infty\), that is, \(q=r=\infty\), then \eqref{GN-Ap} follows from \eqref{pointwise}, since \(\norm{Mf}_{L^\infty(\RR^N)}\leq\norm{f}_{L^\infty(\RR^N)}\) and \(L^\infty(w)=L^\infty(\RR^N)\). If \(p<\infty\), then \eqref{pointwise} and H{\"o}lder's inequality with exponents \(\frac{2q}p\) and \(\frac{2r}p\) give
\[
\int_{\RR^N}\abs{\nabla u}^pw\D x\leq(2\beta_N)^p\int_{\RR^N}(Mu)^{\frac p2}\pt{M(\Delta u)}^{\frac p2}w\D x\leq(2\beta_N)^p\norm{Mu}_{L^q(w)}^{\frac p2}\norm{M(\Delta u)}_{L^r(w)}^{\frac p2},
\]
and \eqref{GN-Ap} follows. If \(w\in A_1\), let \(G_\lambda\) and \(f=\lambda^2u-\Delta u\) be as in the proof of \cref{lem:pointwise}. Fubini's theorem and the symmetry of \(\abs{\nabla G_\lambda}\) give
\[
\int_{\RR^N}\abs{\nabla u}w\D x\leq\int_{\RR^N}\abs f\pt{\abs{\nabla G_\lambda}\ast w}\D x\leq\frac{\beta_N}\lambda\int_{\RR^N}\abs f\,Mw\D x\leq\beta_N[w]_{A_1}\pt{\lambda\norm{u}_{L^1(w)}+\frac1\lambda\norm{\Delta u}_{L^1(w)}},
\]
and we conclude by choosing \(\lambda\) as in Proof \#3.
\end{proof}

\begin{remark}\label{rem:Ap}
For \(w\equiv1\), \cref{gn-thm-Ap} gives \cref{gn-thm-ND} for \(1<q,r\leq\infty\) and for \(p=q=r=1\), with \(\Delta u\) in place of \(D^2u\). Pointwise multiplicative inequalities of this type, with the gradient of the highest order in place of \(\Delta u\), were obtained by Ka{\l}amajska \cite{Kal1994}*{Theorem 3}, who used them to derive Gagliardo--Nirenberg inequalities in weighted Sobolev spaces, and by Maz'ya and Shaposhnikova \cite{MaSh1999}. In \cite{Kal1994}*{Theorem 8}, the derivatives of the highest order may be replaced by a family of homogeneous operators whose characteristic polynomials have no common zero in \(\mathbb C^N\setminus\{0\}\). For \(N\geq2\) this excludes the single operator \(\Delta\), since \(\xi_1^2+\dots+\xi_N^2\) vanishes at \((1,i,0,\dots,0)\). For \(N\geq3\), Maz'ya and Shaposhnikova \cite{MaSh2000} prove \(M(I_zf)\leq C\pt{M(I_\zeta f)}^{\frac z\zeta}(Mf)^{1-\frac z\zeta}\) for the Riesz potentials \(I_z=(-\Delta)^{-\frac z2}\), \(0<z<\zeta<N\). With \(f=-\Delta u\), \(z=1\) and \(\zeta=2\), this is an inequality of the form \eqref{pointwise}, with \((-\Delta)^{1/2}u\) in place of \(\nabla u\). In dimension one, \eqref{pointwise} reads \(\abs{u'(x)}^2\leq4\,Mu(x)\,Mu''(x)\). Maz'ya and Shaposhnikova obtained this inequality with the constant \(8\) in \cite{MaSh1999}*{inequality (4)}, and later the inequality \(\abs{u'(x)}^2\leq\frac83\,M^\diamond u(x)\,M^\diamond u''(x)\), with the sharp constant \(\frac83\), where \(M^\diamond u(x):=\sup_{r>0}\frac1{2r}\abs{\int_{x-r}^{x+r}\operatorname{sgn}(y-x)\,u(y)\D y}\) \cite{MaSh2002}*{inequality (16)}. Since \(M^\diamond u\leq Mu\), this gives \(\abs{u'(x)}^2\leq\frac83\,Mu(x)\,Mu''(x)\). Ka{\l}amajska \cite{Kal2002} gives numerical bounds for the constants \(C_{k,m}\) in the one-dimensional inequalities \(\abs{f^{(k)}}\leq C_{k,m}(Mf)^{1-\frac km}(Mf^{(m)})^{\frac km}\), which are optimal for \(k=1\) and \(m=2\). With the Hessian in place of \(\Delta u\), Maz'ya and Shaposhnikova also obtain the constant \(4(N+1)\) in every dimension, with a modified maximal operator in place of \(Mu\) \cite{MaSh1999}*{Remark 2}. \Cref{gn-thm-Ap} applies to every weight comparable to \(\abs x^\alpha\) with \(-N<\alpha<N(\min\{q,r\}-1)\). In particular, it applies to every homogeneous weight of degree \(\alpha\) in this range that is continuous and positive on the unit sphere. It also applies to the monomial weights \(\abs{x_1}^{A_1}\cdots\abs{x_N}^{A_N}\) with \(-1<A_i<\min\{q,r\}-1\) for every \(i\), since a product \(w(x)=w_1(x_1)\cdots w_N(x_N)\) of one-dimensional \(A_s\) weights belongs to \(A_s\). Indeed, the averages of \(w\) and \(w^{-\frac1{s-1}}\) over a cube \(Q=I_1\times\cdots\times I_N\) are the products of the averages of \(w_i\) and \(w_i^{-\frac1{s-1}}\) over \(I_i\), and every ball \(B\) is contained in a cube \(Q\) with \(\abs Q=c_N\abs B\), so that \([w]_{A_s}\leq c_N^s\prod_i[w_i]_{A_s}\).
\end{remark}

\begin{remark}\label{rem:power}
Let \(1<q,r<\infty\). Since \(\abs x^\alpha\in L^1_{\mathrm{loc}}(\RR^N)\) if and only if \(N_\alpha>0\), the natural range for \eqref{GN-w} with \(w=\abs x^\alpha\) and \(u\in C^\infty_c(\RR^N)\) is \(\alpha>-N\). \Cref{gn-thm-hom} gives \eqref{GN-w} with \(D^2u\) for \(\alpha>r-N\), and \cref{gn-thm-Ap} gives it for \(-N<\alpha<N(\min\{q,r\}-1)\). If \(N\geq2\) and \(r<Nq\) (for instance, if \(q\geq r\)), these two ranges cover all \(\alpha>-N\). If \(N=1\) and \(q\geq r\), they cover all \(\alpha>-1\) except \(\alpha=r-1\).

For power weights, weighted Gagliardo--Nirenberg inequalities with higher derivatives were proved by C.-S. Lin \cite{Lin1986}, extending \cite{CKN1984}. In our setting, the hypotheses of Lin's theorem, as stated in \cite{DuSi2022}*{Theorem 1.2}, reduce to \(N_\alpha>0\) and \(r\neq N\). See also \cite{DuSi2022} for results under Muckenhoupt-type conditions; for power weights, \cite{DuSi2022}*{Theorem 1.6} gives \eqref{GN-w} in the range of \cref{gn-thm-Ap}. The interest of \cref{gn-thm-Ap,gn-thm-hom} is that they follow from the elementary arguments of \cref{sec:p=q=r}, with explicit constants.

Every homogeneous weight that is positive and continuous on the unit sphere is comparable to \(\abs x^\alpha\). For such weights, \eqref{GN-w} follows from the power case, with a constant depending on \(\max_{\abs x=1}w/\min_{\abs x=1}w\).
\end{remark}

\begin{remark}\label{rem:weights-other}
\emph{The measure-adapted operator.} If \(w\) is positive and of class \(C^1\) in an open set \(\Omega\subseteq\RR^N\), then the proof of \cref{gn-thm-RN}, with integration by parts with respect to \(w\D x\), gives, for every \(u\in C^\infty_c(\Omega)\) and with no further assumption on \(w\),
\begin{gather*}
\norm{\nabla u}_{L^p(w)}^2\leq\norm{u}_{L^q(w)}\norm{\abs{L_wu}+(p-2)\abs{D^2u}}_{L^r(w)}\qquad\text{if }p\geq2,\\
\norm{\nabla u}_{L^p(w)}^2\leq\frac{r'}q\norm{u}_{L^q(w)}\norm{L_wu}_{L^r(w)}\qquad\text{if }p\leq2\text{ and }r>1,
\end{gather*}
where for \(p\geq2\) one uses the identity
\[
\operatorname{div}\pt{w\abs{\nabla u}^{p-2}\nabla u}=w\abs{\nabla u}^{p-2}\pt{L_wu+(p-2)D^2u(\nu,\nu)}.
\]
Consider the monomial weights with \(A_i\geq0\), and \(\Omega=\{x:x_i>0\text{ whenever }A_i>0\}\). The boundary terms vanish because \(w=0\) on \(\partial\Omega\), so these inequalities hold for every \(u\in C^\infty_c(\overline\Omega)\), as in \cite{CR-O2013-2}.

\emph{Product weights.} The slicing argument of the proof of \cref{gn-thm-ND} extends to product weights \(w(x)=w_1(x_1)\cdots w_N(x_N)\). If every \(w_i\) satisfies the one-dimensional inequality
\[
\norm{v'}_{L^p(w_i)}^2\leq C_i\norm{v}_{L^q(w_i)}\norm{v''}_{L^r(w_i)},
\]
then the same computation gives
\[
\norm{\partial_iu}_{L^p(w)}^2\leq C_i\norm{u}_{L^q(w)}\norm{\partial_{ii}u}_{L^r(w)}.
\]
For the monomial weights, which degenerate on the sphere and are not covered by \cref{rem:power}, this reduces \eqref{GN-w} to the one-dimensional weights \(\abs t^{A}\), \(A>-1\). For these, the one-dimensional inequality holds if \(A>r-1\) (\cref{gn-thm-hom}) or \(A<\min\{q,r\}-1\) (\cref{gn-thm-Ap}), and in fact, if \(1<q,r<\infty\), for every \(A>-1\) by \cite{Lin1986}. Functions defined only on \(\overline\Omega\) can be treated by first extending them across the hyperplanes \(\{x_i=0\}\) with a reflection of class \(C^{1,1}\).

\emph{Proof \#4.} If \(w\in A_p\) with \(1<p<\infty\), then \(w(x)w(y)\) is an \(A_p\) weight in \(\RR^{2N}\). So the weighted Calder{\'o}n--Zygmund inequality
\[
\norm{D^2v}_{L^p(w\otimes w)}\leq C\norm{\Delta v}_{L^p(w\otimes w)}
\]
holds \cite{CoiFe1974}, and the computation of Proof \#4 gives
\[
\norm{\nabla u}_{L^p(w)}^2\leq C\norm{u}_{L^p(w)}\norm{\Delta u}_{L^p(w)}.
\]
\end{remark}

\begin{remark}\label{rem:general-w}
\Cref{lem:general-a} also applies to a weighted measure \(\mu=w\D x\) on \(\RR^N\), with \cref{gn-thm-hom} giving \eqref{half-lem}, as soon as \(w\D x\) supports a Sobolev inequality. For the power weights \(w=\abs x^\alpha\) with \(\alpha\geq0\), the Caffarelli--Kohn--Nirenberg inequalities \cite{CKN1984} give
\[
\norm{v}_{L^{r^*}(w)}\leq S\norm{\nabla v}_{L^r(w)},\qquad1\leq r<N_\alpha,\quad\frac1{r^*}=\frac1r-\frac1{N_\alpha},
\]
for all Lipschitz functions \(v\) with compact support. Hence, if \(1<r<N_\alpha\), \cref{lem:general-a}, with \(n=N_\alpha\), \(g=\abs{\nabla u}\) and \(h=\abs{D^2u}\), and \cref{gn-thm-hom} give
\[
\norm{\nabla u}_{L^p(w)}\leq C\norm{u}_{L^q(w)}^{1-a}\norm{D^2u}_{L^r(w)}^a,\qquad\frac1p=\frac1{N_\alpha}+a\pt{\frac1r-\frac2{N_\alpha}}+\frac{1-a}q,
\]
for all \(1\leq q\leq\infty\) and \(\frac12\leq a\leq1\), that is, \eqref{gen-GN-N} with \(N_\alpha\) in place of \(N\). Similarly, the monomial weights with \(A_i\geq0\) satisfy a Sobolev inequality with \(N+A_1+\cdots+A_N\) in place of \(N_\alpha\) \cite{CR-O2013-2}, and the same argument applies whenever the case \(a=\frac12\) is available (see \cref{rem:weights-other}).
\end{remark}

\appendix

\section{Proof of \texorpdfstring{\cref{lp-norm-equiv}}{Lemma 1}}

First, observe that by H{\"o}lder's inequality
\[
\abs{v}_{L^p(\RR)}\leq \norm{v}_{L^p(\RR)},
\]
for all \(v\in L^p(\RR)\), so it is enough to find a constant \(C>0\) such that \(\norm{v}_{L^p(\RR)}\leq C\abs{v}_{L^p(\RR)}\) for every \(v\). To this end, it suffices to exhibit \(\varphi\in W^{2,p'}(\RR)\) and a constant \(C>0\) such that
\begin{equation}\label{esti3}
\norm{v}_{L^p(\RR)}\leq C\frac{\displaystyle\int_\RR v(x)\spt{\varphi(x)-\varphi''(x)}\D x}{\norm{\varphi}_{L^{p'}(\RR)}+\norm{\varphi''}_{L^{p'}(\RR)}}.
\end{equation}
Fix \(v\in L^p(\RR)\), let \(f:=\abs{v}^{p-2}v\in L^{p'}(\RR)\), and let \(\varphi\) be the unique solution of the ODE
\begin{equation*} %\label{ode}
\left\{
\begin{aligned}
-\varphi''+\varphi&=f &&\text{in }(-\infty,\infty),\\
\lim_{\abs x\to\infty}\varphi(x)&=0.
\end{aligned}
\right.
\end{equation*}
As in Proof \#3, \(\varphi=\Phi\ast f\), where \(\Phi\) is given by \eqref{defn-phi}. In particular, \(\Phi\in W^{1,1}(\RR)\) with
\[
\norm{\Phi}_{L^1(\RR)}=\norm{\Phi'}_{L^1(\RR)}=1,
\]
and since \(\varphi=\Phi\ast f\) and \(\varphi'=\Phi'\ast f\), we can use Young's inequality for convolutions to obtain
\begin{gather*}
\norm{\varphi}_{L^{p'}(\RR)}\leq \norm{\Phi}_{L^1(\RR)}\norm{f}_{L^{p'}(\RR)}\leq \norm{f}_{L^{p'}(\RR)},\\
\norm{\varphi'}_{L^{p'}(\RR)}\leq \norm{\Phi'}_{L^1(\RR)}\norm{f}_{L^{p'}(\RR)}\leq \norm{f}_{L^{p'}(\RR)},\\
\norm{\varphi''}_{L^{p'}(\RR)}=\norm{\varphi-f}_{L^{p'}(\RR)}\leq 2\norm{f}_{L^{p'}(\RR)}.
\end{gather*}
Therefore, \(\varphi\in W^{2,p'}(\RR)\) and
\[
\frac{\displaystyle\int_\RR v(x)\spt{\varphi(x)-\varphi''(x)}\D x}{\norm{\varphi}_{L^{p'}(\RR)}+\norm{\varphi''}_{L^{p'}(\RR)}}=\frac{\norm{v}_{L^p(\RR)}^p}{\norm{\varphi}_{L^{p'}(\RR)}+\norm{\varphi''}_{L^{p'}(\RR)}}\geq \frac{\norm{v}_{L^p(\RR)}^p}{3\norm{\abs{v}^{p-1}}_{L^{p'}(\RR)}}=\frac13 \norm{v}_{L^p(\RR)},
\]
and the result follows with \(C=3\).

\qed

\section*{Acknowledgments}

This work started around 2010, while I was a Ph.D. student at Rutgers University under the supervision of Haim Brezis, and I have returned to it on and off since then. I am deeply grateful to the late Professor Brezis for suggesting the study of this problem and for encouraging me to find these proofs. I also thank H. Wang for very interesting discussions during that period, from which Proofs \#3 and \#4 originated.

\section*{Declarations}

\noindent\textbf{Competing interests.} The author declares no competing interests.

\noindent\textbf{Data availability.} Data sharing is not applicable to this article, as no datasets were generated or analysed during the current study.

\bibliography{references-B}

@Article{BrMi2018,
  author  = {Brezis, Haim and Mironescu, Petru},
  title   = {Gagliardo--{N}irenberg inequalities and non-inequalities: the full story},
  journal = {Ann. Inst. H. Poincar\'{e} Anal. Non Lin\'{e}aire},
  year    = {2018},
  volume  = {35},
  number  = {5},
  pages   = {1355--1376},
  doi     = {10.1016/j.anihpc.2017.11.007},
}

@Article{Cas2016-2,
  author   = {Castro, Hern{\'a}n},
  title    = {Hardy-{S}obolev-type inequalities with monomial weights},
  journal  = {Ann. Mat. Pura Appl. (4)},
  fjournal = {Annali di Matematica Pura ed Applicata. Series IV},
  year     = {2017},
  volume   = {196},
  number   = {2},
  pages    = {579--598},
  issn     = {0373-3114},
  doi      = {10.1007/s10231-016-0587-2},
}

@Article{Cas2021,
  author   = {Castro, Hern{\'a}n},
  title    = {Extremals for {H}ardy-{S}obolev type inequalities with monomial weights},
  journal  = {J. Math. Anal. Appl.},
  fjournal = {Journal of Mathematical Analysis and Applications},
  year     = {2021},
  volume   = {494},
  number   = {2},
  pages    = {Paper No. 124645},
  issn     = {0022-247X},
  doi      = {10.1016/j.jmaa.2020.124645},
}

@Article{CaZy1952,
  author  = {Calder\'{o}n, A. P. and Zygmund, A.},
  title   = {On the existence of certain singular integrals},
  journal = {Acta Math.},
  year    = {1952},
  volume  = {88},
  pages   = {85--139},
  doi     = {10.1007/BF02392130},
}

@Article{CFR2020,
  author  = {Ciraolo, Giulio and Figalli, Alessio and Roncoroni, Alberto},
  journal = {Geom. Funct. Anal.},
  title   = {Symmetry results for critical anisotropic {$p$}-{L}aplacian equations in convex cones},
  year    = {2020},
  number  = {3},
  pages   = {770--803},
  volume  = {30},
  doi     = {10.1007/s00039-020-00535-3},
}

@Article{CiMa2022,
  author   = {Cialdea, A. and Maz'ya, V.},
  title    = {A survey of functional and {$L^p$}-dissipativity theory},
  journal  = {Bull. Math. Sci.},
  fjournal = {Bulletin of Mathematical Sciences},
  year     = {2022},
  volume   = {12},
  number   = {2},
  pages    = {Paper No. 2230003},
  issn     = {1664-3607},
  doi      = {10.1142/S1664360722300031},
}

@Article{CKN1984,
  author   = {Caffarelli, Luis A. and Kohn, Robert V. and Nirenberg, Louis},
  title    = {First order interpolation inequalities with weights},
  journal  = {Compositio Math.},
  fjournal = {Compositio Mathematica},
  year     = {1984},
  volume   = {53},
  number   = {3},
  pages    = {259--275},
  issn     = {0010-437X},
}

@Article{CoDu2003,
  author  = {Coulhon, Thierry and Duong, Xuan Thinh},
  title   = {Riesz transform and related inequalities on noncompact {R}iemannian manifolds},
  journal = {Comm. Pure Appl. Math.},
  year    = {2003},
  volume  = {56},
  number  = {12},
  pages   = {1728--1751},
}

@Article{CoiFe1974,
  author   = {Coifman, R. R. and Fefferman, C.},
  title    = {Weighted norm inequalities for maximal functions and singular integrals},
  journal  = {Studia Math.},
  fjournal = {Polska Akademia Nauk. Instytut Matematyczny. Studia Mathematica},
  year     = {1974},
  volume   = {51},
  pages    = {241--250},
  issn     = {0039-3223,1730-6337},
  doi      = {10.4064/sm-51-3-241-250},
}

@Article{CR-O2013-2,
  author   = {Cabr{\'e}, Xavier and Ros-Oton, Xavier},
  title    = {Sobolev and isoperimetric inequalities with monomial weights},
  journal  = {J. Differential Equations},
  fjournal = {Journal of Differential Equations},
  year     = {2013},
  volume   = {255},
  number   = {11},
  pages    = {4312--4336},
  issn     = {0022-0396},
  doi      = {10.1016/j.jde.2013.08.010},
}

@Article{CRS2016,
  author  = {Cabr{\'e}, Xavier and Ros-Oton, Xavier and Serra, Joaquim},
  journal = {J. Eur. Math. Soc. (JEMS)},
  title   = {Sharp isoperimetric inequalities via the {ABP} method},
  year    = {2016},
  number  = {12},
  pages   = {2971--2998},
  volume  = {18},
  doi     = {10.4171/JEMS/659},
}

@Article{DuSi2022,
  author  = {Duarte, Rodrigo and Drumond Silva, Jorge},
  title   = {Weighted {G}agliardo--{N}irenberg interpolation inequalities},
  journal = {J. Funct. Anal.},
  year    = {2023},
  volume  = {285},
  number  = {5},
  pages   = {Paper No. 110009},
  doi     = {10.1016/j.jfa.2023.110009},
}

@Book{Eva2010,
  author    = {Evans, Lawrence C.},
  title     = {Partial differential equations},
  edition   = {2},
  series    = {Graduate Studies in Mathematics},
  volume    = {19},
  publisher = {American Mathematical Society},
  address   = {Providence, RI},
  year      = {2010},
  doi       = {10.1090/gsm/019},
}

@Article{FiFoRoSo2021,
  author  = {Fiorenza, Alberto and Formica, Maria Rosaria and Roskovec, Tom\'{a}\v{s} G. and Soudsk\'{y}, Filip},
  title   = {Detailed proof of classical {G}agliardo--{N}irenberg interpolation inequality with historical remarks},
  journal = {Z. Anal. Anwend.},
  year    = {2021},
  volume  = {40},
  number  = {2},
  pages   = {217--236},
  doi     = {10.4171/ZAA/1681},
}

@Article{Gag59,
  author  = {Gagliardo, Emilio},
  title   = {Ulteriori propriet\`{a} di alcune classi di funzioni in pi\`{u} variabili},
  journal = {Ricerche Mat.},
  year    = {1959},
  volume  = {8},
  pages   = {24--51},
}

@Book{GT01,
  author    = {Gilbarg, David and Trudinger, Neil S.},
  title     = {Elliptic partial differential equations of second order},
  series    = {Classics in Mathematics},
  publisher = {Springer-Verlag},
  address   = {Berlin},
  year      = {2001},
  note      = {Reprint of the 1998 edition},
  isbn      = {3-540-41160-7},
  doi       = {10.1007/978-3-642-61798-0},
}

@Article{GuPi2019,
  author  = {G\"{u}neysu, Batu and Pigola, Stefano},
  title   = {{$L^p$}-interpolation inequalities and global {S}obolev regularity results},
  journal = {Ann. Mat. Pura Appl. (4)},
  year    = {2019},
  volume  = {198},
  number  = {1},
  pages   = {83--96},
  note    = {With an appendix by Ognjen Milatovic},
  doi     = {10.1007/s10231-018-0763-7},
}

@Book{HL52,
  author    = {Hardy, G. H. and Littlewood, J. E. and P{\'o}lya, G.},
  title     = {Inequalities},
  edition   = {2},
  publisher = {Cambridge, at the University Press},
  year      = {1952},
}

@Article{Kal1994,
  author  = {Ka{\l}amajska, Agnieszka},
  title   = {Pointwise multiplicative inequalities and {N}irenberg type estimates in weighted {S}obolev spaces},
  journal = {Studia Math.},
  year    = {1994},
  volume  = {108},
  number  = {3},
  pages   = {275--290},
  doi     = {10.4064/sm-108-3-275-290},
}

@Article{Kal2002,
  author  = {Ka{\l}amajska, Agnieszka},
  title   = {On constants in the pointwise one-dimensional multiplicative inequalities},
  journal = {Comment. Math. (Prace Mat.)},
  year    = {2002},
  volume  = {42},
  number  = {1},
  pages   = {37--51},
}

@Article{KaPe2012,
  author  = {Ka{\l}amajska, Agnieszka and Peszek, Jan},
  title   = {On some nonlinear extensions of the {G}agliardo--{N}irenberg inequality with applications to nonlinear eigenvalue problems},
  journal = {Asymptot. Anal.},
  year    = {2012},
  volume  = {77},
  number  = {3-4},
  pages   = {169--196},
  doi     = {10.3233/ASY-2011-1079},
}

@Article{KaPi2006,
  author  = {Ka{\l}amajska, Agnieszka and Pietruska-Pa{\l}uba, Katarzyna},
  title   = {Gagliardo--{N}irenberg inequalities in weighted {O}rlicz spaces},
  journal = {Studia Math.},
  year    = {2006},
  volume  = {173},
  number  = {1},
  pages   = {49--71},
  doi     = {10.4064/sm173-1-4},
}

@Article{Kol1939,
  author  = {Kolmogorov, A. N.},
  title   = {On inequalities between the upper bounds of the successive derivatives of an arbitrary function on an infinite interval},
  journal = {Uchen. Zap. Moskov. Gos. Univ. Matematika},
  year    = {1939},
  volume  = {30},
  pages   = {3--16},
  note    = {In Russian; English translation in Amer. Math. Soc. Transl. (1) 2 (1962), 233--243},
}

@Book{KwZe1992,
  author    = {Kwong, Man Kam and Zettl, Anton},
  title     = {Norm inequalities for derivatives and differences},
  series    = {Lecture Notes in Mathematics},
  volume    = {1536},
  publisher = {Springer-Verlag, Berlin},
  year      = {1992},
  doi       = {10.1007/BFb0090864},
}

@Article{Lan1913,
  author  = {Landau, Edmund},
  title   = {Einige {U}ngleichungen f\"{u}r zweimal differentiierbare {F}unktionen},
  journal = {Proc. London Math. Soc. (2)},
  year    = {1914},
  volume  = {13},
  pages   = {43--49},
  doi     = {10.1112/plms/s2-13.1.43},
}

@Article{LeRoSo2025,
  author  = {Le{\'s}nik, Karol and Roskovec, Tom{\'a}{\v{s}} G. and Soudsk{\'y}, Filip},
  title   = {{G}agliardo--{N}irenberg inequality via a new pointwise estimate},
  journal = {J. Funct. Anal.},
  year    = {2025},
  volume  = {289},
  number  = {7},
  pages   = {Paper No. 110996},
  doi     = {10.1016/j.jfa.2025.110996},
}

@Article{Lin1986,
  author  = {Lin, Chang-Shou},
  title   = {Interpolation inequalities with weights},
  journal = {Comm. Partial Differential Equations},
  year    = {1986},
  volume  = {11},
  number  = {14},
  pages   = {1515--1538},
  doi     = {10.1080/03605308608820473},
}

@Article{MaSh1999,
  author  = {Maz'ya, Vladimir and Shaposhnikova, Tatyana},
  title   = {On pointwise interpolation inequalities for derivatives},
  journal = {Math. Bohem.},
  year    = {1999},
  volume  = {124},
  number  = {2-3},
  pages   = {131--148},
  doi     = {10.21136/MB.1999.126252},
}

@InCollection{MaSh2000,
  author    = {Maz'ya, Vladimir and Shaposhnikova, Tatyana},
  title     = {Pointwise interpolation inequalities for {R}iesz and {B}essel potentials},
  booktitle = {Analytical and computational methods in scattering and applied mathematics ({N}ewark, {DE}, 1998)},
  series    = {Chapman \& Hall/CRC Res. Notes Math.},
  volume    = {417},
  publisher = {Chapman \& Hall/CRC, Boca Raton, FL},
  year      = {2000},
  pages     = {217--229},
}

@Article{MaSh2002,
  author  = {Maz'ya, Vladimir and Shaposhnikova, Tatyana},
  title   = {Sharp pointwise interpolation inequalities for derivatives},
  journal = {Funct. Anal. Appl.},
  year    = {2002},
  volume  = {36},
  number  = {1},
  pages   = {30--48},
  doi     = {10.1023/A:1014478100799},
}

@Article{MeSp2008,
  author  = {Metafune, Giorgio and Spina, Chiara},
  title   = {An integration by parts formula in {S}obolev spaces},
  journal = {Mediterr. J. Math.},
  year    = {2008},
  volume  = {5},
  number  = {3},
  pages   = {357--369},
  doi     = {10.1007/s00009-008-0155-0},
}

@Article{Muck1972,
  author   = {Muckenhoupt, Benjamin},
  title    = {Weighted norm inequalities for the {H}ardy maximal function},
  journal  = {Trans. Amer. Math. Soc.},
  fjournal = {Transactions of the American Mathematical Society},
  year     = {1972},
  volume   = {165},
  pages    = {207--226},
  issn     = {0002-9947},
  doi      = {10.1090/S0002-9947-1972-0293384-6},
}

@Article{Nir59,
  author  = {Nirenberg, L.},
  title   = {On elliptic partial differential equations},
  journal = {Ann. Scuola Norm. Sup. Pisa (3)},
  year    = {1959},
  volume  = {13},
  pages   = {115--162},
}

@Misc{Ouh2024,
  author = {Ouhabaz, El Maati},
  title  = {A multiplicative inequality of {R}iesz transform type on general {R}iemannian manifolds},
  note   = {Preprint, arXiv:2406.01097},
  year   = {2024},
}

@Article{Ste1983,
  author  = {Stein, Elias M.},
  title   = {Some results in harmonic analysis in {$\mathbb{R}^n$}, for {$n\to\infty$}},
  journal = {Bull. Amer. Math. Soc. (N.S.)},
  year    = {1983},
  volume  = {9},
  number  = {1},
  pages   = {71--73},
  doi     = {10.1090/S0273-0979-1983-15157-1},
}

@Book{Ste70,
  author    = {Stein, Elias M.},
  title     = {Singular integrals and differentiability properties of functions},
  series    = {Princeton Mathematical Series, No. 30},
  publisher = {Princeton University Press, Princeton, N.J.},
  year      = {1970},
  doi       = {10.1515/9781400883882},
}

\end{document}